\documentclass[leqno]{amsart}
\usepackage{egastyle}
\usepackage{mymacros}
\usepackage[all,tips,dvips]{xy}
\usepackage{url}
\usepackage{bm}
\usepackage{verbatim}

\def\T{\bT}
\def\G{\bG}
\def\A{\bA}

\DeclareMathOperator{\val}{val}
\DeclareMathOperator{\trop}{trop}
\DeclareMathOperator{\Trop}{Trop}
\DeclareMathOperator{\relint}{relint}
\DeclareMathOperator{\Int}{Int}
\DeclareMathOperator{\conv}{conv}
\DeclareMathOperator{\inn}{in}
\newcommand{\an}{{\operatorname{an}}}
\newcommand{\sTor}{\mathscr{T}\mathrm{or}}

\newcommand{\s}[1]{{}^{#1}}
\newcommand{\p}{{}'}

\CompileMatrices

\title{Lifting non-proper tropical intersections}

\author{Brian Osserman}
\email{osserman@math.ucdavis.edu}
\address{Department of Mathematics, One Shields Avenue, University of
  California, Davis, CA 95616}
\author{Joseph Rabinoff}
\email{rabinoff@math.harvard.edu}
\address{Department of Mathematics, One Oxford Street, Harvard University,
  Cambridge, MA 02138} 

\thanks{The first author was partially supported by NSA grant H98230-11-1-0159, and the second author was supported by an NSF postdoctoral research fellowship.}

\begin{document}

\maketitle

\begin{abstract}  \vspace{-30pt}
  We prove that if $X,X'$ are closed subschemes of a torus $\T$ over a
  non-Archimedean field $K$, of
  complementary codimension and with finite intersection, then the 
  stable tropical intersection along a (possibly
  positive-dimensional, possibly unbounded) connected component
  $C$ of $\Trop(X)\cap\Trop(X')$ lifts to algebraic
  intersection points, with multiplicities.  This theorem requires potentially passing
  to a suitable toric variety $X(\Delta)$ and its associated
  extended tropicalization $N_\R(\Delta)$; the
  algebraic intersection points lifting the stable tropical intersection
  will have tropicalization somewhere in the closure of $C$ in $N_\R(\Delta)$.  
  The proof involves a result on continuity of intersection numbers in
  the context of non-Archimedean analytic spaces.
\end{abstract}

\section{Introduction}
Let $K$ be a field equipped
with a nontrivial%
\footnote{In the introduction we assume that the valuation is nontrivial
  for simplicity; we will prove the main theorems in the trivially valued
  case as well.} 
non-Archimedean valuation $\val:K\to\R\cup\{\infty\}$,
and suppose that $K$ is complete or algebraically closed. Let
$\T\cong\G_m^n$ be a finite-rank split torus over $K$ with coordinate
functions $x_1,\ldots,x_n$.  The
\emph{tropicalization} map is the function
$\trop:|\T|\to\R^n$ given by 
$\trop(\xi) = (\val(x_1(\xi)),\ldots,\val(x_n(\xi)))$, where $|X|$
denotes the set of closed points of a scheme $X$.  Given a closed
subscheme $X\subseteq\T$, the \emph{tropicalization} of $X$ is the closure
(with respect to the Euclidean topology) of the set
$\trop(|X|)$ in $\R^n$, and is denoted $\Trop(X)$.  This is a subset
which can be endowed with the structure of a weighted polyhedral complex
of the same dimension as $X$.  In particular, it is a combinatorial
object, a ``shadow'' of $X$ which is often much easier to analyze than $X$
itself.  It is therefore important that one can recover information about
$X$ from its tropicalization.

An example of this idea is to relate the intersection of $X$ with a second
closed subscheme $X'\subseteq\T$ to the intersections of their
tropicalizations.  One might hope that
$\Trop(X\cap X') = \Trop(X)\cap\Trop(X')$, but this is not generally the
case.  For example, let $K$ be the field of Puiseux series over $\C$ with
uniformizer $t$.  The curves
$X=\{ x + y = 1 \}$ and $X' = \{ tx + y = 1 \}$ do not meet in $\G_m^2$,
but $\Trop(X)\cap\Trop(X')$ is the ray 
$\R_{\geq 0}\cdot (1,0)\subset\R^2$.  This example is
``degenerate'' in the sense that $\Trop(X)$ does not intersect $\Trop(X')$
transversely; generically the intersection of two
one-dimensional polyhedral complexes in $\R^2$ is a finite set of
points.  This is in fact the only obstruction: assuming $X,X'$ pure
dimensional, if $\Trop(X)$ meets
$\Trop(X')$ in the expected codimension at a point $v\in\R^n$, then
$v\in\Trop(X\cap X')$.  This was proved by Osserman and
Payne, who in fact prove much
more: they show that in a suitable sense, the tropicalization of the
intersection cycle $X\cdot X'$ is equal to the stable tropical
intersection $\Trop(X)\cdot\Trop(X')$, still under the hypothesis that
$\Trop(X)$ meets $\Trop(X')$ in the expected codimension; see
\cite[Theorem~1.1, Corollary~5.1.2]{osserman_payne:lifting}.
In particular, if $\codim(X)+\codim(X')= \dim(\T)$ and
$\Trop(X)\cap\Trop(X')$ is a finite set of points, then 
$\Trop(X)\cdot\Trop(X')$ is a weighted sum of points of $\R^n$; these
points then lift, with multiplicities, to points of $X\cdot X'$.
Hence in this case one can compute local intersection numbers via
tropicalization. 

This paper will be concerned with the case when
$\codim(X)+\codim(X')=\dim(\T)$, but when the intersection
$\Trop(X)\cap\Trop(X')$ may
have higher-dimensional connected components.
The stable tropical intersection $\Trop(X)\cdot\Trop(X')$ is still a
well-defined finite set of points contained in $\Trop(X)\cap\Trop(X')$,
obtained by translating $\Trop(X)$ by a generic vector $v$ and then taking
the limit as $v\to 0$, but it is no longer the case that 
$\Trop(X\cdot X') = \Trop(X)\cdot\Trop(X')$.  Indeed, in the above example
of $X = \{x+y=1\}$ and $X'=\{tx+y=1\}$, the stable tropical intersection
is the point $(0,0)$ with multiplicity $1$, but $X\cap X'=\emptyset$.
This illustrates the need to compactify the situation in the direction of
the ray $\R_{\geq 0}\cdot(1,0)$.  Let us view $X,X'$ as curves in 
$\A^1\times\G_m$, and extend the tropicalization map to 
a map $\trop: |\A^1\times\G_m|\to (\R\cup\{\infty\})\times\R$ in the
obvious way.  Then $X\cap X'$ is the reduced point $(0,1)$, and
$\Trop(X\cap X') = \{(\infty,0)\}$ is contained in the closure of $\Trop(X)\cap\Trop(X')$
in $(\R\cup\{\infty\})\times\R$.  It is not a coincidence that
the multiplicity of the point $(0,1)\in X\cap X'$ coincides with the
multiplicity of $(0,0)\in\Trop(X)\cdot\Trop(X')$: we have lost the
ability to pinpoint the exact location of the point $\Trop(X\cap X')$
beyond saying that it lies in the closure of $\Trop(X)\cap\Trop(X')$, but
we are still able to recover its multiplicity using the stable tropical
intersection.

In order to carry out this strategy in general, we need to make precise
the notion of ``compactifying in the directions where
the tropicalization is infinite'':
we say that an integral pointed fan $\Delta$ is a 
\emph{compactifying fan} for a polyhedral complex $\Pi$
provided that the recession cone of each cell of $\Pi$ is a union of cones
in $\Delta$.  The setup for the main theorem is then as follows.
Let $X_1,\ldots,X_m\subseteq\T$ be pure-dimensional closed subschemes with 
$\sum_{i=1}^m \codim(X_i)=\dim(\T)$,
and let $C\subseteq\bigcap_{i=1}^m \Trop(X_i)$ be a connected component.  
Let $\Pi$ be the polyhedral complex underlying $C$ (with respect to
some choice of polyhedral complex structures on the $\Trop(X_i)$), and let
$\Delta$ be a compactifying fan for $\Pi$.
Let $M$ be the lattice of
characters of $\T$ and let $N$ be its dual lattice, so
the $\Trop(X_i)$ naturally live in $N_\R = N\tensor_\Z\R$.
We partially compactify the torus with the toric variety $X(\Delta)$,
which contains $\T$ as a dense open subscheme.  The
extended tropicalization is a topological space
$N_\R(\Delta)$ which canonically contains $N_\R$ as a dense open subset,
and which is equipped with a map 
$\trop: |X(\Delta)|\to N_\R(\Delta)$ extending 
$\trop: |\T|\to N_\R$; see~\parref{par:extended.trop}.
Let $\bar C$ be the closure of $C$ in $N_\R(\Delta)$;
this is a compact set since $\Delta$ is a compactifying fan for $\Pi$
(Remark~\ref{rem:cpct.fans.cpct}). 
For an isolated point $\xi\in\bigcap_{i=1}^m\bar X_i$ we let
$i_K(\xi,\bar X_1 \cdots \bar X_m; X(\Delta))$ denote the multiplicity of $\xi$
in the intersection class $\bar X_1 \cdots \bar X_m$, and for 
$v\in\bigcap_{i=1}^m \Trop(X_i)$ we
let $i(v,\prod_{i=1}^m \Trop(X_i))$ denote the multiplicity of $v$ in the
stable tropical intersection $\Trop(X_1) \cdots \Trop(X_m)$.

\begin{thm*}
  If $X(\Delta)$ is smooth, and if there are only finitely many points of
  $|\bar X_1\cap\cdots\cap\bar X_m|$ mapping to $\bar C$ under
  $\trop$, then
  \[ \sum_{\substack{\xi\in|\bigcap_{i=1}^m\bar X_i|\\\trop(\xi)\in\bar C}}
  i_K\big(\xi,\, \bar X_1\cdots\bar X_m;\, X(\Delta)\big) = 
  \sum_{v\in C} i\big(v,\, \Trop(X_1)\cdots\Trop(X_m)\big).
  \]
\end{thm*}

See Theorem~\ref{thm:main-multiple}.
This can be seen as a lifting theorem for points in the stable tropical 
intersection $\Trop(X_1) \cdots \Trop(X_m)$, with the provisos that
we may have to do some compactification of the situation first, and that
the tropicalizations of the
points of the algebraic intersection $\bar X_1 \cdots \bar X_m$ corresponding
to a point $v$ of the stable tropical intersection are only
confined to the closure of the connected component of
$\bigcap_{i=1}^m\Trop(X_i)$ containing $v$.

The finiteness assumption on $\bigcap_{i=1}^m \bar X_i$ is also necessary
in this generality --- we will provide some conditions under which it is
automatically satisfied. In particular, when $C$ is \emph{bounded}, 
the compactifying fan is unnecessary, and we have:

\begin{cor*}
  Suppose 
  that $C$ is 
  bounded.  Then there are only finitely many points of 
  $|X_1\cap\cdots\cap X_m|$ mapping to $C$ under
  $\trop$, and
  \[ \sum_{\substack{\xi\in|\bigcap_{i=1}^m X_i|\\\trop(\xi)\in C}}
  i_K\big(\xi,\, X_1\cdots X_m;\, \T\big) = 
  \sum_{v\in C} i\big(v,\, \Trop(X_1)\cdots\Trop(X_m)\big).
  \]
\end{cor*}

The proof of the main Theorem proceeds as follows.  Assume for simplicity
that $K$ is both complete and algebraically closed.  Let
$X,X'\subseteq\T$ be pure-dimensional closed subschemes with 
$\codim(X)+\codim(X')=\dim(\T)$,
and let $C$ be a connected component of $\Trop(X)\cap\Trop(X')$.  
Assume for the moment that $C$ is bounded.  Let 
$v\in N$ be a generic cocharacter, regarded as a homomorphism
$v:\G_m\to\T$.  Then $(\Trop(X)+\epsilon\cdot v)\cap\Trop(X')$ is a finite
set for small enough $\epsilon$, and the stable tropical intersection is
equal to $\lim_{\epsilon\to 0}(\Trop(X)+\epsilon\cdot v)\cap\Trop(X')$;
this can be seen as a ``continuity of local tropical intersection numbers''.
For $t\in K^\times$ the tropicalization of $v(t)\cdot X$ is equal to
$\Trop(X) - \val(t)\cdot v$, so for small nonzero values of $\val(t)$ we can
apply Osserman-Payne's tropical lifting theorem to 
$(v(t)\cdot X)\cap X'$.  Hence what we want to prove is a theorem on
continuity of \emph{local} algebraic intersection numbers that applies to the family 
$\fY = \{(v(t)\cdot X)\cdot X'\}_{\val(t)\in[-\epsilon,\epsilon]}$.

There are two problems with proving this continuity of local intersection
numbers, both of which have the same solution.  The first is that the base
of the family $\fY$ is the set 
$\cS_\epsilon(K) = \{t\in K^\times~:~\val(t)\in[-\epsilon,\epsilon]\}$,
which is not algebraic but an \emph{analytic} annulus in $\G_m$.  The
second is that we only want to count 
intersection multiplicities in a neighborhood of $C$ --- more precisely, if $P$
is a polytope containing $C$ in its interior and disjoint from the other
components of $\Trop(X)\cap\Trop(X')$, then for every
$t\in\cS_\epsilon(K)$ we only want to count intersection multiplicities of
points in $\cU_P(K) = \trop\inv(P)$, which is again an analytic subset of
$\T$.  Therefore we will prove that dimension-zero intersection numbers of
analytic spaces are constant in flat families over an analytic base.
This is one of the main ideas of the paper; the other idea, orthogonal to
this one, is the precise
compactification procedure described above, which is necessary when $C$ is
unbounded.

\medskip
\noindent\textbf{Outline of the paper.}\,~\,Many of the technical
difficulties in this paper revolve around the need 
to pass to a compactifying toric variety when our connected component $C$
of $\Trop(X)\cap\Trop(X')$ is unbounded.  As such, section~\ref{sec:fans}
is devoted to introducing compactifying and compatible fans $\Delta$, and
studying the behavior of the closure operation for polyhedra in
$N_\R(\Delta)$.  The main result is
Proposition~\ref{prop:compactifying.fans}, which says in particular that 
for a suitable fan $\Delta$, the extended tropicalization of the intersection 
of the closures of $X$ and $X'$ in $X(\Delta)$ is contained in the closure of 
$\Trop(X)\cap\Trop(X')$, and that the same can be achieved for individual
connected components of the intersection.
This is quite important in the statement of the
Theorem above, since we want to sum over all closed points
$\xi$ of $\bar X\cap\bar X\p$ with $\trop(\xi)\in\bar C$, and is also
vital in section~\ref{sec:moving}.

In section~\ref{sec:moving} we prove a version of the tropical moving
lemma: the stable tropical intersection $\Trop(X)\cdot\Trop(X')$ is
defined locally by translating $\Trop(X)$ by a small amount in the
direction of a generic displacement vector $v$, and in
Lemma~\ref{lem:moving.lemma} we make these conditions precise.  The main
point of section~\ref{sec:moving}, however, is to show that for 
$v\in N$ satisfying the tropical moving lemma, the corresponding family 
$\{(v(t)\cdot\bar X)\cap\bar X\p\cap\cU_P\}_{t\in\cS_\epsilon}$
of analytic subspaces of $\cU_P$, where $P$ is a polyhedral
neighborhood of $C$, is \emph{proper} over $\cS_\epsilon$.  See
Proposition~\ref{prop:conn.comp}.  We therefore
give a brief discussion of the analytic notion of properness
in section~\ref{sec:moving}, which we conclude with the
very useful tropical criterion for properness of a family of analytic
subspaces of a toric variety (Proposition~\ref{prop:properness}).

In section~\ref{sec:int.nums} we define local intersection multiplicities
of dimension-zero intersections of analytic spaces in a smooth ambient
space, using a slight modification of Serre's definition.  These analytic
intersection numbers coincide with the algebraic ones in the case of
analytifications of closed subschemes
(Proposition~\ref{prop:alg.an.int.nums}).  The main result
(Proposition~\ref{prop:cont_nums}) is the
continuity of analytic intersection numbers mentioned above: if
$\cX,\cX'$ are analytic spaces, flat over a connected base $\cS$, inside a smooth
analytic space $\cZ$, such that $\cX\cap\cX'$ is finite over $\cS$, then
the total intersection multiplicities on any two fibers are equal.

In section~\ref{sec:mainthm} we prove the main theorem
(Theorem~\ref{thm:mainthm}) and its corollaries, combining the results of
sections~\ref{sec:moving} and~\ref{sec:int.nums}.  We
also treat the case of intersecting more than two subschemes of $\T$ by
reducing to intersection with the diagonal.  We conclude by giving a
detailed worked example in section~\ref{sec:the.example}.

\section{Analytifications and tropicalizations}
We will use the following general notation throughout the paper.  If
$P$ is a subspace of a topological space $X$, its interior (resp.\
closure) in $X$ will be denoted $P^\circ$ (resp.\ $\bar P$).  If
$f:X\to Y$ is a map (of sets, schemes, analytic spaces, etc.) the fiber
over $y\in Y$ will be denoted $X_y = f\inv(y)$.

By a \emph{cone} in a Euclidean space we will always mean a polyhedral cone.

\paragraph[Non-Archimedean fields]
We fix a non-Archimedean field $K$, i.e.\ a field 
equipped with a non-Archimedean valuation $\val: K\to\R\cup\{\infty\}$.  We will
assume throughout that $K$ is complete or algebraically closed,
and except in \parref{par:trop}, \parref{par:extended.trop}, and
section~\ref{sec:mainthm}, we assume further 
that $\val$ is nontrivial and that $K$ is complete with respect to
$\val$, in order to be able to work with analytic spaces over $K$.
Let $|\cdot| = \exp(-\val(\cdot))$ be the corresponding absolute value
and let $G = \val(\bar K\s\times)\subseteq\R$ be the saturation of the value
group of $K$.  

By a \emph{valued field extension} of $K$ we mean a non-Archimedean field
$K'$ equipped with an embedding $K\inject K'$ which respects the
valuations.

\paragraph[Analytic spaces] \label{par:analytic.spaces}
Assume that $K$ is complete and nontrivially valued.%
\footnote{For convenience, in this paper we only work with analytic spaces
  over nontrivially valued fields, although Berkovich's theory is valid in
  the trivially valued case.}
In this paper, by an \emph{analytic space} we mean a separated (i.e.\
the underlying topological space is Hausdorff), good,
strictly $K$-analytic space in the sense
of~\cite{berkovich:etalecohomology}.  In particular, all $K$-affinoid
algebras and $K$-affinoid spaces are assumed to be strictly $K$-affinoid.
We will generally use calligraphic
letters to refer to analytic spaces.  For a $K$-affinoid algebra $A$,
its Berkovich spectrum $\sM(A)$ is an analytic space whose underlying
topological space is the set of bounded multiplicative semi-norms 
$\|\cdot\|:A\to\R_{\geq 0}$, equipped with topology of pointwise
convergence.  An affinoid space is compact.
If $\cX$ is an analytic space, $|\cX|$ will denote the set
of classical ``rigid'' points of $\cX$; this definition is local on $\cX$, and if
$\cX = \sM(A)$ is affinoid, then $|\cX|$ is naturally identified with the
set of maximal ideals of $A$.  The subset $|\cX|$ is everywhere dense
in $\cX$ by~\cite[Proposition~2.1.15]{berkovich:analytic_geometry}.  We
also let $\cX(\bar K) = \varinjlim_{K'} \cX(K')$, where
$\cX(K') = \Hom_K(\sM(K'),\cX)$ and the union runs over
all finite extensions $K'$ of $K$.  There is a natural surjective map 
$\cX(\bar K)\surject|\cX|$.  

For a point $x$ of an analytic space $\cX$, we let $\sH(x)$ denote the
completed residue field at $x$.  This is a complete valued field extension of $K$
which plays the role of the residue field at a point of a scheme.  In
particular, if $\cY\to\cX$ is a morphism, then the set-theoretic fiber
$\cY_x$ is naturally an $\sH(x)$-analytic space.  The
point $x$ is rigid if and only if $[\sH(x):K] < \infty$, in 
which case we will generally use the notation $K(x) = \sH(x)$.

Let $\cX$ be an analytic space.  An \emph{analytic domain} in $\cX$ is,
roughly, a subset $\cY$ which naturally inherits the structure of analytic
space from $\cX$.  These play the role of the open subschemes of a scheme;
in particular, any open subset of $\cX$ is an analytic domain.  An
analytic domain need not be open, however; for example, an \emph{affinoid domain} in
$\cX$ is an analytic domain which is also an 
affinoid space (which is compact, hence closed).  A 
\emph{Zariski-closed subspace} of $\cX$ is an analytic 
space $\cY\inject\cX$ which is locally defined by the vanishing of some number of
analytic functions on $\cX$.  The set underlying $\cY$ is closed in $\cX$.

For any separated, finite-type $K$-scheme $X$ we let $X^\an$ denote the
analytification of $X$.  This analytic space comes equipped with a map of
ringed spaces $X^\an\to X$ which identifies the set $|X|$ of closed points
(resp.\ the set $X(\bar K)$ of geometric points) with
$|X^\an|$ (resp.\ $X^\an(\bar K)$).  If for $x\in|X|$ we let $K(x)$ denote
the residue field at $x$, then $K(x)$ is identified with the completed
residue field of the associated point $x\in|X^\an|$.
The analytification functor respects all fiber products and complete valued
extensions of the ground field.
In the case that $X = \Spec(A)$ is
affine, we will identify the topological space underlying $X^\an$ with
the space of all multiplicative semi-norms
$\|\cdot\|:A\to\R\cup\{\infty\}$ extending the absolute value on $K$.

If $X$ is a $K$-scheme (resp.\ a $K$-analytic space) and $K'$ is a field
extension (resp.\ complete valued field extension) of $K$, we let $X_{K'}$ denote
the base change to $K'$. 

\paragraph[Tropicalization]\label{par:trop}
Here we assume that $K$ is a complete or algebraically closed,
possibly trivially-valued non-Archimedean field.
Let $M\cong\Z^n$ be a finitely generated free
abelian group and $N = \Hom_\Z(M,\Z)$ its dual.  
For any subgroup $\Gamma\subseteq\R$ we let $M_\Gamma = M\tensor_\Z \Gamma$ and  
$N_\Gamma = N\tensor_\Gamma \Z = \Hom_\Z(M,\Gamma)$.
Let $\T = \Spec(K[M])$ be the torus with character lattice $M$.
Given a closed subscheme $X \subseteq \T$ and a point $v \in N_\R$,
the \emph{initial degeneration} $\inn_v(X)$ is a canonically defined 
scheme over the residue field of $K$, of finite type if 
$v \in N_G \subseteq N_\R$. The \emph{tropicalization}
of $X$ is the subset $\Trop(X) \subseteq N_\R$ of all $v$ such that $\inn_v(X)$
is nonempty. If $K'/K$ is a complete or algebraically closed valued field
extension then  
$\Trop(X_{K'}) = \Trop(X)$.  The set $\Trop(X)$ can be enriched with the
structure of a polyhedral complex (which is in general non-canonical) 
with the property that if $v,v' \in N_G$ lie in the interior of the
same cell, then $\inn_v(X_{\bar K}) \cong \inn_{v'}(X_{\bar K})$. This polyhedral complex
has positive integer weights canonically assigned to each facet, defined
as follows: let $P\subset\Trop(X)$ be a facet, let 
$v\in\relint(P)$, and let $K'$ be an algebraically closed
valued field
extension of $K$ with value group $G'$ such that $v\in N_{G'}$.  The
\emph{tropical multiplicity} $m(P)$ of
$P$ is defined to be the sum of the multiplicities of the irreducible
components of $\inn_v(X_{K'})$.  This is independent of the choice of $K'$ 
by~\cite[Remark~A.5]{osserman_payne:lifting},
or~\cite[\S4.18]{bpr:nonarch_trop} in the complete case (see
also~\cite[Lemma~4.19]{bpr:nonarch_trop}).  
The weights are insensitive to algebraically closed valued field
extensions. 
See for instance~\cite[\S2]{osserman_payne:lifting}
for a more detailed survey of the above. 

If $X = V(f)$ is the hypersurface defined by a Laurent polynomial $f\in K[M]$
then we write $\Trop(f) = \Trop(X)$; the set $\Trop(f)$ is equipped with a
canonical weighted polyhedral complex structure.  See for
instance~\cite[\S8]{jdr:trop_ps}. 

For $u\in M$, let $x^u\in K[M]$ denote the corresponding character.
The \emph{tropicalization map} $\trop: |\T|\to N_\R$ is the map defined by
$\angles{u,\trop(\xi)} = -\val(x^u(\xi))$, where
$\angles{\cdot,\cdot}:M_\R\times N_\R\to\R$ is the canonical pairing.
We also denote the composition $\T(\bar K)\to|\T|\to N_\R$ by $\trop$.
Note that this definition only makes sense when $K$ is complete or
algebraically closed, as $\val(x^u(\xi))$ is not in general well-defined
if $K$ is neither.  If $K$ is complete and nontrivially valued, we define
$\trop:\T^\an\to N_\R$ by $\angles{u,\trop(\|\cdot\|)} = \log(\|x^u\|)$;
this is a continuous, proper surjection which
is compatible with $\trop:|\T|\to N_\R$ under the identification
$|\T|=|\T^\an|$.  

Let $X\subseteq\T$ be a closed subscheme.
If $K$ is nontrivially valued then $\Trop(X)$ is the closure of 
$\trop(|X|)$ in $N_\R$, and if in addition $K$ is complete then
$\Trop(X) = \trop(X^\an)$.  If $K$ is trivially valued then
$\trop(|X|) = \{0\}$ or is empty.

\paragraph[Extended tropicalization]\label{par:extended.trop}
We continue to assume that $K$ is a complete or algebraically closed,
possibly trivially-valued non-Archimedean field.
If $\sigma$ is an integral cone in $N_\R$, we let $X(\sigma)$ denote the
affine toric variety $\Spec(K[\sigma^\vee\cap M])$, where 
\[ \sigma^\vee = \{ u\in M_\R~:~\angles{u,v}\leq 0
\text{ for all } v\in\sigma \}. \]
Likewise for an integral fan $\Delta$ in $N_\R$ we let 
$X(\Delta)$ be the toric variety obtained by gluing the affine toric
varieties $X(\sigma)$ for $\sigma\in\Delta$.

Let $\sigma$ be an integral cone in $N_\R$.  We define
$N_\R(\sigma) = \Hom_{\R_{\geq0}}(\sigma^\vee,\R\cup\{-\infty\})$, the
set of homomorphisms of additive monoids with an action of 
$\R_{\geq 0}$.  We equip $N_\R(\sigma)$ with the topology of pointwise
convergence.  The
tropicalization map extends to a map
$\trop:|X(\sigma)|\to N_\R(\sigma)$, again using the formula 
$\angles{u,\trop(\xi)} = -\val(x^u(\xi))$.  If $K$ is complete and
nontrivially valued then we define
$\trop:X(\sigma)^\an\to N_\R(\sigma)$ by 
$\angles{u,\trop(\|\cdot\|)} = \log(\|x^u\|)$; as above this is a
continuous, proper surjection which is compatible with
$\trop:|X(\sigma)|\to N_\R(\sigma)$ under the identification
$|X(\sigma)| = |X(\sigma)^\an|$.
If $\Delta$ is an integral fan in $N_\R$ we set
$N_\R(\Delta) = \bigcup_{\sigma\in\Delta} N_\R(\sigma)$; the tropicalization maps 
$\trop:|X(\sigma)|\to N_\R(\sigma)$ 
(resp.\ $\trop:X(\sigma)^\an\to N_\R(\sigma)$ in the complete nontrivially
valued case) glue to give a map
$\trop:|X(\Delta)|\to N_\R(\Delta)$
(resp.\ a continuous, proper surjection
$\trop:X(\Delta)^\an\to N_\R(\Delta)$).
As above we also use $\trop$ to denote the composite map 
$X(\Delta)(\bar K)\to |X(\Delta)| \to N_\R(\Delta)$.  

There is a natural decomposition  
$N_\R(\Delta) = \Djunion_{\sigma\in\Delta} N_\R/\spn(\sigma)$, which
respects the decomposition of $X(\Delta)$ into torus orbits.  We will make
this identification implicitly throughout the paper. 
If $X(\Delta)=X(\sigma)$ is an affine toric variety, a monoid
homomorphism $v:\sigma^\vee\to\R\cup\{-\infty\}$ is in the stratum
$N_\R/\spn(\tau)$ if and only if $v\inv(\R)=\tau^\perp\cap\sigma^\vee$.
For a cone $\sigma\subseteq N_\R$ we let $\pi_\sigma$ denote the quotient
map $N_\R\to N_\R/\spn(\sigma)$.  
We will use the following explicit description of the
topology on $N_\R(\sigma)$:

\begin{lem}\label{lem:extended-closure} 
  Let $\sigma \subset N_\R$ be a pointed cone.
  A sequence $v_1,v_2,\ldots\in N_\R$ converges to the point
  $\bar v \in N_\R/\spn(\tau)\subseteq N_\R(\sigma)$ for some
  $\tau\prec\sigma$ if and only if both of the following hold:
  \begin{enumerate}
  \item $\angles{u,v_i}\to \angles{u,\bar v}$ as $i\to\infty$
    for all $u\in\sigma^\vee\cap\tau^\perp$ (equivalently, 
    $\pi_\tau(v_i)\to\bar v$ as $i\to\infty$), and
  \item $\angles{u,v_i}\to -\infty$ as $i\to\infty$ 
    for all $u\in\sigma^\vee\setminus\tau^\perp$.
  \end{enumerate}
\end{lem}

\pf Since $N_\R(\sigma)$ is equipped with the topology of pointwise
convergence, this follows immediately from the fact that for
$u\in\sigma^\vee$ we 
have $\angles{u,\bar v} \neq -\infty$ if and only if $u\in\tau^\perp$
(note that since $\sigma$ is pointed, $\sigma^\vee$ spans $M_\R$).\qed

If $X\subseteq X(\Delta)$ is a closed subscheme, its 
\emph{extended tropicalization} $\Trop(X,\Delta)\subset N_\R(\Delta)$ can be defined by 
tropicalizing each torus orbit separately. If the valuation on $K$ is
nontrivial then $\Trop(X,\Delta)$ is the closure of $\trop(|X|)$ in
$N_\R(\Delta)$, and if in addition $K$ is complete then
$\Trop(X,\Delta) = \trop(X^\an)$.  
See~\cite{payne:analytification,jdr:trop_ps} for details on extended
tropicalizations.

\section{Compatible and compactifying fans}\label{sec:fans}
If $\cP$ is any finite collection of polyhedra, its \emph{support} is the
closed subset $|\cP| = \bigcup_{p\in\cP} P$.
In this section we develop the related notions of compatible and
compactifying fans for $\cP$.  Roughly,
if $\Delta$ is compatible with $\cP$ then the closure of 
$|\cP|$ is easy to calculate in $N_\R(\Delta)$, and if $\Delta$ is a
compactifying fan then the closure of $|\cP|$ in
$N_\R(\Delta)$ is compact --- i.e., $N_\R(\Delta)$ compactifies $N_\R$ in
the directions in which $|\cP|$ is infinite.  This will be important when
$|\cP|$ is a connected component of the intersection of tropicalizations.

The \emph{recession cone} of a polyhedron $P\subseteq N_\R$ is defined to be
the set 
\[ \rho(P) = \{w\in N_\R~:~ v+w\in P\text{ for all } v\in P \}. \]
If $P$ is cut out by conditions $\left<u_i,v\right> \leq c_i$ for
$u_1,\dots,u_m \in M_{\R}$ and $c_1,\dots,c_m \in \R$, then
$\rho(P)$ is given explicitly by $\left<u_i,v\right> \leq 0$ for
$i=1,\dots,m$.

\begin{defn}
  Let $\cP$ be a finite collection of polyhedra in $N_\R$ and let 
  $\Delta$ be a pointed fan.
  \begin{enumerate}
  \item The fan $\Delta$ is said to be \emph{compatible with $\cP$}
    provided that, for all $P\in\cP$ and all cones
    $\sigma\in\Delta$, either $\sigma \subseteq \rho(P)$  
    or $\relint(\sigma) \cap \rho(P)=\emptyset$.
  \item The fan $\Delta$ is said to be a \emph{compactifying fan for $\cP$}
    provided that, for all $P\in\cP$, the recession cone $\rho(P)$
    is a union of cones in $\Delta$.
  \end{enumerate}
\end{defn}

The reason that we will generally require our fans to be pointed is
due to the fact that if $\Delta$ is a pointed fan in $N_\R$, then $N_\R$ is
canonically identified with the open subspace $N_\R(\{0\})$ of
$N_\R(\Delta)$. 

Following are some basic properties of compatible and compactifying fans,
which are easily checked directly from the definitions.

\begin{prop} \label{prop:basic-fan-props} 
Let $\cP$ be a finite collection of polyhedra in $N_\R$.

\begin{enumerate}
\item A compactifying fan for $\cP$ is compatible with $\cP$.
\item A subfan of a fan compatible with $\cP$ is compatible with $\cP$.
\item A refinement of a fan compatible with $\cP$ is compatible with
  $\cP$, and a refinement of a compactifying fan for $\cP$ is a
  compactifying fan for $\cP$.
\item If a fan is compatible with $\cP$, it is compatible with any subset
of $\cP$. A compactifying fan for $\cP$ is a compactifying fan for any
subset of $\cP$.
\item\label{item:subset-complex} Suppose that $\cP$ is a subset of the cells of a polyhedral complex
$\Pi$, and $\cP$ contains all the maximal cells of $\Pi$ (equivalently,
$\cP$ and $\Pi$ have the same support). Then a fan is compatible with $\cP$
if and only if it is compatible with $\Pi$, and a fan is a compactifying fan
for $\cP$ if and only if it is a compactifying fan for $\Pi$.
\item Let $\cP'$ be a second finite collection of polyhedra in $N_\R$.
  A fan compatible with both $\cP$ and $\cP'$
  is compatible with $\cP \cap \cP'$, and 
  a compactifying fan for both $\cP$ and $\cP'$ is 
  a compactifying fan for $\cP\cap\cP'$.
\end{enumerate}
\end{prop}

Here the notation $\cP \cap \cP'$ means the set of intersections of pairs
of polyhedra in $\cP$ and $\cP'$. If $\cP$ and $\cP'$
are the sets of cells of polyhedral complexes $\Pi$ and $\Pi'$, then 
$\cP \cap \cP'$ is not generally equal to $\Pi \cap \Pi'$, as it does not
have to contain every face of every polyhedron. However, according 
to Proposition~\ref{prop:basic-fan-props}(\ref{item:subset-complex})
above, a fan is compatible with  
$\cP \cap \cP'$ if and only it is compatible with $\Pi \cap \Pi'$, and a
fan is a compactifying fan for 
$\cP \cap \cP'$ if and only it is a compactifying fan for $\Pi \cap \Pi'$.

\begin{rem} \label{rem:cpct.fans.cpct}
  Let $\cP$ be a finite collection of polyhedra in $N_\R$ and let $\Delta$
  be a compactifying fan for $\cP$.  We claim that the closure of $|\cP|$
  in $N_\R(\Delta)$ is compact.  To prove this we may assume that
  $\cP=\{P\}$ consists of a single polyhedron, and by
  Lemma~\ref{lem:delta.polyhedra} below we may even assume that
  $\rho(P)\in\Delta$.  The closure of $P$ in $N_\R(\rho(P))$ is compact
  by~\cite[\S3]{jdr:trop_ps}, so the claim follows since $N_\R(\rho(P))$
  is a subspace of $N_\R(\Delta)$.
\end{rem}

\begin{defn}\label{defn:cont-family} 
Let $a,b\in\R$ with $a\leq b$ and let $V$ be a finite-dimensional real
vector space.
A \emph{continuous family of polyhedra} in $V$, parameterized by $[a,b]$,
is a function $\cP$ from $[a,b]$ to the set of all polyhedra in $V$, given
by an equation of the form
$$\cP(t) = \bigcap_{i=1}^m \{v \in N_{\R}:\left<u_i,v\right>\leq f_i(t)\},$$
where $u_i \in V^*$  for $i=1,\dots,m$, and $f_i(t)$ a continuous
real-valued function on $[a,b]$.
\end{defn}

Note that in the above definition, we allow $V=(0)$, in which case each
$u_i$ is necessarily $0$, and each $\sP(t)$ is either empty or $V$
according to whether or not all the $f_i(t)$ are nonnegative. In addition,
we allow $a=b$, in which case $\sP$ is just
a polyhedron.  Note also that if $\sP,\sP'$ are continuous families of
polyhedra in $V$ parameterized by $[a,b]$ then $t\mapsto\sP(t)\cap\sP'(t)$
is one also.

For the convenience of the reader we include proofs  of
the following lemmas on polyhedra, which are undoubtedly
well known. 
The first lemma roughly says that if $P$ is a polyhedron, then we have
$\lim_{t\to 0}(tP) = \rho(P)$.

\begin{lem} \label{lem:converge.to.rhoP}
  Let $V$ be a finite-dimensional real vector space and let
  $P\subseteq V$ be a polyhedron.  For $t\in[0,1]$ define
  \[ \sP(t) =
  \begin{cases}
    tP &\quad t \in (0,1] \\
    \rho(P) &\quad t = 0.
  \end{cases}\]
  Then $\sP$ is a continuous family of polyhedra.
\end{lem}

\pf Suppose that $P$ is defined by the inequalities
$\left<u_i,v\right> \leq c_i$ for some $u_1,\dots,u_m \in M_{\R}$ and 
$c_1,\dots,c_m \in \R$.  Then $tP$ is defined by
$\angles{u_i,v}\leq tc_i$ for $i=1,\ldots,m$, so the lemma follows because 
$\rho(P)$ is given by $\angles{u_i,v} \leq 0$ for $i=1,\dots,m$.\qed

\begin{lem} \label{lem:polyhedral-projs} Given a finite-dimensional real
vector space $V$ and a continuous family of polyhedra $\sP$ in $V$, the 
image of $\sP$ under projection to any quotient space $W$ of $V$ is a 
continuous family of polyhedra.
\end{lem}

Here the projection is taken one $t$ at a time, in the obvious way.

\pf Since every projection can be factored as a composition of
projections with $1$-dimensional kernels, it is enough to consider this case.
Accordingly, let $W$ be a quotient of $V$, with the kernel of 
$V \twoheadrightarrow W$ being $1$-dimensional.
Choose a basis $x_1,\dots,x_n$ of $V^*$, and write
$$u_i=\sum_{j=1}^n a_{i,j} x_j$$
for each $i$. We may further suppose that we have chosen the $x_i$ so that
the kernel of the given projection is precisely the intersection of the
kernels of $x_2,\dots,x_n$. Thus, $x_2,\dots,x_n$ gives a basis for $W^*$.
Without loss of generality, we may reorder the $u_i$ so that
$a_{1,1},\dots,a_{p,1}=0$, $a_{p+1,1},\dots,a_{q,1}>0$, and
$a_{q+1,1},\dots,a_{m,1}<0$. Dividing through the $u_i$ and $f_i$ for $i>p$ 
by $a_{i,1}$, we have that the inequalities defining $\sP$ can be
rewritten as follows: for $i=1,\dots,p$, we have
$\left<u'_i,v\right> \leq f_i(t)$,
where $u'_i=u_i$;
for $i=p+1,\dots,q$, we have
$\left<x_1,v\right> \leq \left<u'_i,v\right> + f_i(t)$, 
where $u'_i =x_1-u_i$;
and for $i=q+1,\dots,m$, we have
$\left<x_1,v\right> \geq \left<u'_i,v\right> + f_i(t)$, 
where $u'_i = x_1-u_i$.
Noting that each $u'_i$ is now well defined on $W$, we see that the image
of $\sP$ in $W$ is described by the inequalities
$\left<u'_i,v\right> \leq f_i(t)$
for $i=1,\dots,p$, and
$\left<u'_i-u'_j,v\right> \leq f_j(t)-f_i(t)$
for each $i=q+1,\dots,m$ and $j=p+1,\dots,q$.
We thus conclude the desired statement.\qed

 The following two corollaries of the
lemma will be useful to us. Setting $a=b$ in Lemma \ref{lem:polyhedral-projs}
we have:

\begin{cor}\label{cor:polyhedral-projs} Let $P$ be a polyhedron in
$N_{\R}$. Then $\pi_\sigma(P)$ is a polyhedron, and in particular is closed.
\end{cor}

\medskip

On the other hand, considering projection to the $0$-space we immediately
conclude:

\begin{cor}\label{cor:polyhedral-family} The set of $t$ for which a 
continuous family of polyhedra is nonempty is closed in $[a,b]$.
\end{cor}

\smallskip
Our main lemma is then the following.

\begin{lem} \label{lem:closure.P}
  Let $P$ be a polyhedron in $N_\R$ and let $\Delta$ be a pointed fan.
  If $\bar P$ is the closure of $P$ in $N_\R(\Delta)$, then
  \[\bar P = 
  \Djunion_{\substack{\sigma\in\Delta\\\relint(\sigma)\cap\rho(P)\neq \emptyset}}
  \pi_\sigma(P). \]
\end{lem}

\pf Since $N_{\R}(\sigma)$ embeds as an open subset of $N_{\R}(\Delta)$
for each $\sigma$ in $\Delta$, and $N_{\R}(\tau)$ embeds as an open
subset of $N_{\R}(\sigma)$ for each face $\tau$ of $\sigma$, to analyze
$\bar P$ it suffices to consider the stratum of $\bar P$ lying in
$N_{\R}/\spn(\sigma)$ for each $\sigma \in \Delta$.

Let $\rho = \rho(P)$.
First suppose $\relint(\sigma) \cap \rho \neq \emptyset$, and let
$\bar v\in\pi_\sigma(P)$. Choose 
$v\in P$ such that $\pi_\sigma(v) = \bar v$, and let
$w\in\relint(\sigma)\cap\rho$. 
Then $aw\in\rho$ for all $a\in\R_{>0}$, so $v + aw\in P$ by the definition
of $\rho$.  But $v + aw\to\bar v$ as $a\to\infty$ by
Lemma~\ref{lem:extended-closure}, so we have $\bar v\in\bar P$.
Hence we obtain one containment.

Next, suppose that $\bar v \in N_{\R}(\sigma)/\spn(\sigma)$ is in
$\bar P$.  Then according to Lemma~\ref{lem:extended-closure} there exists
a sequence $v_1,v_2,\dots \in P$ such that 
$\lim_{i \to \infty} \pi_{\sigma}(v_i) = \bar v$,
and for all $u\in\sigma^\vee\setminus\sigma^{\perp}$
we have $\lim_{i \to \infty} \left<u,v_i\right> = -\infty$. In particular,
we see that $\bar v$ is in the closure of $\pi_{\sigma}(P)$, hence in 
$\pi_{\sigma}(P)$ by Corollary~\ref{cor:polyhedral-projs}. It is therefore
enough to show that $\relint(\sigma) \cap \rho \neq \emptyset$.

Choose generators $u_1,\dots,u_m$ for $\sigma^{\vee}$.
For $\delta\geq 0$, denote by $\sigma_{\delta}$ the polyhedron cut out by the 
conditions
$\left<u_j,v\right> \leq -1$ for $u_j \not\in \sigma^{\perp}$, and
$\left<u_j,v\right> \leq \delta$ for $u_j \in \sigma^{\perp}$. Then 
$\sigma_0 \subseteq \relint(\sigma)$.
Fix $\delta > 0$ and choose $\epsilon > 0$ such that 
$\epsilon \left<u_j,\bar v\right> < \delta$ for all $j$ such that
$u_j\in\sigma^\perp$. 
For all $i\gg 0$ we have $\left<u_j,\epsilon v_i\right> \leq -1$ when
$u_j \not\in \sigma^{\perp}$ (since 
$\angles{u_j,\epsilon v_i}\to -\infty$) and 
$\angles{u_j,\epsilon v_i}\leq\delta$ when
$u_j\in\sigma^\perp$ (since 
$\angles{u_j,\epsilon v_i}\to\epsilon\angles{u_j,\bar v}$).
We thus see that for fixed $\delta$ and
$\epsilon$ sufficiently small, we have 
$(\epsilon P) \cap \sigma_{\delta} \neq \emptyset$.
Still holding $\delta$ fixed, by Lemma~\ref{lem:converge.to.rhoP} we see
that 
\[ \epsilon\mapsto
\begin{cases}
  (\epsilon P) \cap \sigma_{\delta} &\quad \epsilon \in(0,1] \\
  \rho\cap\sigma_{\delta} &\quad \epsilon = 0
\end{cases}\]
forms a continuous family of polyhedra,
so by Corollary \ref{cor:polyhedral-family} we conclude that
$\rho \cap \sigma_{\delta} \neq \emptyset$. But now
letting $\delta$ vary, we have that
$\rho \cap \sigma_{\delta}$ also forms a continuous family of polyhedra,
so $\rho \cap \sigma_0 \neq \emptyset$, and $\rho$ meets the relative
interior of $\sigma$, as desired.\qed

\begin{lem} \label{lem:closure.int}
  Let $\cP,\cP'$ be finite collections of polyhedra and let $\Delta$ be a 
  pointed fan in $N_\R$.  If $\Delta$ is
  compatible with either $\cP$ or $\cP'$ then
  \[ \bar{|\cP|\cap|\cP'|} = \bar{|\cP|}\cap\bar{|\cP'|}, \]
  where all closures are taken in $N_\R(\Delta)$.
\end{lem}

\pf We assume without loss of generality that $\Delta$ is compatible with
$\cP$.  Let $P\in\cP$ and $P'\in\cP'$.
It suffices to show that
$\bar{P\cap P'} = \bar P\cap\bar P\p$.
First we claim that $\pi_\sigma(P)\cap\pi_\sigma(P') = \pi_\sigma(P\cap P')$
for all $\sigma\in\Delta$ such that
$\relint(\sigma)\cap\rho(P)\cap\rho(P')\neq\emptyset$; note that this
condition is equivalent to $\sigma\subseteq\rho(P)$ and
$\relint(\sigma)\cap\rho(P')\neq\emptyset$.  It is obvious that
$\pi_\sigma(P)\cap\pi_\sigma(P') \supset \pi_\sigma(P\cap P')$, so let
$\bar v\in\pi_\sigma(P)\cap\pi_\sigma(P')$.  Choose
$v\in P$ and $v'\in P'$ such that
$\pi_\sigma(v) = \pi_\sigma(v') = \bar v$, and choose
$w'\in\relint(\sigma)\cap\rho(P')$.  For any $a\in\R$ we have
$v'+aw'-v\in\spn(\sigma)$, so for $a\gg 0$ we have
$v' +aw' - v\in\sigma\subseteq\rho(P)$.  Choose such an $a$, and set
$w = v' + aw' - v$. Then
$v + w = v' + aw'$; since $v + w\in P$ and $v'+aw'\in P'$, this shows that
\[ \bar v = \pi_\sigma(v+w) = \pi_\sigma(v' + aw') \in
\pi_\sigma(P\cap P'). \]

By Lemma~\ref{lem:closure.P} as applied to $P$ and $P'$, we have
\begin{align*} \bar P\cap\bar P\p
& = \Djunion_{\substack{\sigma\in\Delta\\\sigma\subseteq\rho(P)}} \pi_\sigma(P)
\cap \Djunion_{\substack{\sigma\in\Delta\\\relint(\sigma)\cap\rho(P')\neq\emptyset}}
\pi_\sigma(P')
= \Djunion_{\substack{\sigma\in\Delta\\
    \relint(\sigma)\cap\rho(P)\cap\rho(P')\neq\emptyset}}
\big(\pi_\sigma(P)\cap\pi_\sigma(P')\big) \\
& = \Djunion_{\substack{\sigma\in\Delta\\\relint(\sigma)\cap\rho(P\cap P')\neq\emptyset}}
\big(\pi_\sigma(P\cap P')\big),
\end{align*}
where the last equality follows from the above and the fact that for
$P \cap P' \neq \emptyset$, we have $\rho(P \cap P')=\rho(P) \cap \rho(P')$.
Applying Lemma~\ref{lem:closure.P} to $P \cap P'$, this last expression
is precisely $\bar{P \cap P'}$.\qed

\smallskip
Applying Lemma \ref{lem:closure.int} twice, we obtain the following.

\begin{cor}\label{cor:compat-intersect}
  Let $\cP,\cP',\cQ$ be finite collections of polyhedra, and let $\Delta$ be 
  a pointed fan in $N_\R$. Suppose that $\Delta$ is compatible 
  with $\cQ$ and with either $\cP \cap \cQ$ or $\cP' \cap \cQ$. Then 
  \[ \bar{|\cP|}\cap\bar{|\cP'|}\cap \bar{|\cQ|} = 
 \bar{|\cP|\cap|\cP'|\cap|\cQ|},
\]
  where all closures are taken in $N_\R(\Delta)$.
\end{cor}

\smallskip
Now we apply Corollary~\ref{cor:compat-intersect} to tropicalizations of 
subschemes.  Assume that $K$ is complete and nontrivially valued.

\begin{prop} \label{prop:compactifying.fans}
  \begin{enumerate}
  \item Let $X$ be a closed subscheme of $\T$ and let $\Delta$ be an
    integral pointed fan in $N_\R$.  Let $\bar X$ be the closure
    of $X$ in $X(\Delta)$.  Then $\Trop(\bar X,\Delta)$ is the closure of
    $\Trop(X)$ in $N_\R(\Delta)$.

  \item Let $X,X'$ be closed subschemes of $\T$, let $\cP$ be a 
    finite collection of polyhedra in $N_\R$, and let $\Delta$ be a
    fan compatible with $\cP$ and with either $\Trop(X)\cap \cP$ or 
    $\Trop(X')\cap \cP$.  Then 
    \[ \Trop(\bar X,\Delta)\cap\Trop(\bar X\p,\Delta) \cap \bar{|\cP|} =
    \bar{\Trop(X)\cap\Trop(X')\cap|\cP|} \] 
    in $N_\R(\Delta)$.
  \end{enumerate}
\end{prop}

\pf The first part is~\cite[Lemma~3.1.1]{osserman_payne:lifting}, and the
second part follows immediately from the first part together with
Corollary~\ref{cor:compat-intersect}.
\qed

\begin{rem} \label{rem:cpct.fans.exist}
  If $\cP$ is a finite collection of polyhedra, then there always exists a
  compactifying fan $\Delta$ for $\cP$. Indeed, given $P_i \in \cP$, let
  $\Delta_i$ be a complete fan containing $\rho(P_i)$ (see for
  instance~\cite{rohrer:completions_fans}). Let $\Delta$ be a common
  pointed refinement of all the 
  $\Delta_i$. Then according to Proposition~\ref{prop:basic-fan-props}, 
  $\Delta$ is a compactifying fan for $\cP$. Although 
  $\Delta$ is complete, we may pass to a compactifying fan with minimal
  support by letting $\Delta'$ be the subfan of $\Delta$ consisting of all
  cones contained in $\rho(P)$ for some $P \in \cP$, and it is clear that
  $\Delta'$ is still a compactifying fan for $\cP$. If $\cP$ consists of
  integral polyhedra, then we may choose the $\Delta_i$ and hence $\Delta$
  and $\Delta'$ to be integral as well.

  For the specific case of tropicalizations, we may also proceed as follows.
  If $X = V(f)$ is the hypersurface defined by a nonzero Laurent polynomial 
  $f\in K[M]$,
  then any pointed refinement of the normal fan to the Newton polytope of
  $f$ is a (complete) compactifying fan for $\Trop(X)$;
  see~\cite[\S12]{jdr:trop_ps}.   In general, one appeals to 
  the \emph{tropical basis theorem}, which states that there exist 
  generators $f_1,\ldots,f_r$ of the ideal defining $X$ such that
  $\Trop(X) = \bigcap_{i=1}^r \Trop(f_i)$.
  See~\cite[\S2.5]{maclagan_sturmfels:book}.
  Any fan which simultaneously refines a compactifying fan for each $V(f_i)$
  is a compactifying fan for $\Trop(X)$ by
  Proposition~\ref{prop:basic-fan-props}. As above, such a fan will be 
  complete, but we may always pass to a suitable subfan.
\end{rem}

\section{The moving lemma} \label{sec:moving}
We begin this section by proving a tropical moving lemma, which roughly
says that if $X,X'\subseteq\T$ are closed subschemes with
$\codim(X)+\codim(X')=\dim(\T)$, and if $\Trop(X)\cap\Trop(X')$ is not a
finite set of points, then for any connected component $C$ of 
$\Trop(X) \cap \Trop(X')$ and generic $v\in N$, there exists a
small $\epsilon>0$, and a neighborhood $C'$ of $C$, such that for all
$t\in[-\epsilon,0)\cup(0,\epsilon]$,
the set $(\Trop(X)+t v)\cap\Trop(X')\cap C'$ is finite, and 
furthermore that for all $t \in [-\epsilon,\epsilon]$, the 
intersection of the closures of $(\Trop(X)+t v)$, $\Trop(X')$, and $C'$ 
is precisely the closure of $(\Trop(X)+t v) \cap \Trop(X')$.

The main point of this section is to give an analytic
counterpart to this deformation, in the following sense.
Let $C$ be a connected component of $\Trop(X)\cap\Trop(X')$, and assume
for simplicity that $C$ is bounded.  Let $P$ be a polytope in $N_\R$
containing $C$ in its interior and such that 
$\Trop(X)\cap\Trop(X')\cap P = C$.  We will express the family
$\{(\Trop(X)+t v)\cap\Trop(X')\cap P\}_{t\in[-\epsilon,\epsilon]}$ as the
tropicalization of a natural family $\cY$ of analytic subspaces of $\T^\an$
parameterized by an analytic annulus $\cS$, which we can then study with
algebraic and analytic methods.  The main result of this section is that
$\cY\to\cS$ is \emph{proper}.

Much of the technical difficulty in this section is in treating the case
when $C$ is not bounded.  This requires quite precise control over the
relationships between the various polyhedra and fans which enter the
picture.

\paragraph[The tropical moving lemma]
Let $P$ be an integral $G$-affine polyhedron in $N_\R$, so
$P = \bigcap_{i=1}^r \{v\in N_\R~:~\angles{u_i,v}\leq a_i\}$ for
some $u_1,\ldots,u_r\in M$ and $a_1,\ldots,a_r\in G$.  
As in~\cite[\S12]{jdr:trop_ps}, we define a \emph{thickening} of $P$ to be
a polyhedron of the form
\[ P' = \bigcap_{i=1}^r \{v\in N_\R~:~\angles{u_i,v}\leq a_i+\epsilon
\} \]
for some $\epsilon > 0$ in $G$.  Note that $\rho(P')=\rho(P)$
and that $P$ is contained in the interior $(P\p)^\circ$ of $P'$.
If $\cP$ is a finite collection
of integral $G$-affine polyhedra, a \emph{thickening} of $\cP$ is a
collection of (integral $G$-affine) polyhedra of the form 
$\cP' = \{P'~:~P\in\cP\}$, where $P'$
denotes a thickening of $P$.  

\begin{rem} \label{rem:thickening.interior}
  Let $P$ be a pointed integral $G$-affine polyhedron and let $P'$
  be a thickening of $P$.  Let $\sigma = \rho(P) = \rho(P')$.  Then the closure
  $\bar P$ of $P$ in $N_\R(\sigma)$ is contained in the interior of 
  $\bar P\p$; see Lemma~\ref{lem:closure.P}
  and~\cite[Remark~3.4]{payne:analytification}. 
  More generally, if $\cP$ is a finite collection of integral $G$-affine
  polyhedra with recession cones contained in a pointed fan $\Delta$, and
  if $\cP'$ is a thickening of $\cP$, then the closure 
  $\bar{|\cP|}$ of $|\cP|$ in $N_\R(\Delta)$ is contained in the interior
  of $\bar{|\cP'|}$.
\end{rem}

\begin{defn}
  Let $\Delta$ be an integral pointed fan and 
  let $\cP$ be a finite collection of integral $G$-affine polyhedra in $N_\R$. A 
  \emph{refinement of $\cP$} is a finite collection of integral $G$-affine
  polyhedra $\cP'$ such
  that every polyhedron of $\cP'$ is contained in some polyhedron of $\cP$,
  and every polyhedron of $\cP$ is a union of polyhedra in $\cP'$. A
  \emph{$\Delta$-decomposition of $\cP$} is a
  refinement $\cP'$ of $\cP$ such that $\rho(P)\in\Delta$ for all $P\in\cP'$.
  A \emph{$\Delta$-thickening of $\cP$} is a thickening of a
  $\Delta$-decomposition of $\cP$.
\end{defn}

If $\cP'$ is a refinement of $\cP$ then $|\cP'|=|\cP|$.  If $\cP'$ is a
$\Delta$-thickening of $\cP$ then $\bar{|\cP|}\subseteq\bar{|\cP'|}{}^\circ$
by Remark~\ref{rem:thickening.interior}.

\begin{lem} \label{lem:delta.polyhedra}
  Let $\cP$ be a finite collection of integral $G$-affine polyhedra
  and let $\Delta$ be an integral
  compactifying fan for $\cP$. Then there exists a $\Delta$-decomposition 
  $\cP'$ of $\cP$.  If further $\cP''$ is a finite collection of polyhedra  
  such that $\Delta$ is compatible with $\cP \cap \cP''$, and $\cP'$ is any
  $\Delta$-decomposition of $\cP$, then $\Delta$ is compatible with 
  $\cP' \cap \cP''$. 
\end{lem}

\pf It suffices to prove the first part of the lemma when $\cP = \{P\}$ 
is a polyhedron such that $\rho(P)$ is a union of cones in $\Delta$. 
First suppose that $P$ is pointed, and
let $P_1$ be the convex hull of the vertices of $P$. 
By~\cite[\S3]{jdr:trop_ps} we have $P = P_1 + \rho(P)$,
so $P = \bigcup_{\sigma\in\Delta,\sigma\subseteq\rho(P)} (P_1+\sigma)$.
Hence it is enough to note that if $F$ is an integral $G$-affine polytope
and $\sigma$ is an integral cone then $F+\sigma$ is an integral $G$-affine
polyhedron with recession cone $\sigma$.

Now suppose that $P$ is not pointed.  Let $W'\subseteq\rho(P)$ be the largest
linear space contained in $\rho(P)$ and let $W$ be a complementary integral
subspace in $N_\R$.  Then $P\cap W$ is a pointed polyhedron with recession
cone $\rho(P)\cap W$, so if $P_1$ is the convex hull of the vertices of
$P\cap W$ then $P\cap W = P_1 + \rho(P)\cap W$.  Hence 
$P = P_1 + \rho(P)$, so the proof proceeds as above.

For the second half of the lemma, given $P' \in \cP'$, $P'' \in \cP''$, 
and $\sigma \in \Delta$, suppose that 
$\relint(\sigma) \cap \rho(P' \cap P'') \neq \emptyset$. Then
in particular $P' \cap P'' \neq \emptyset$, and 
$\rho(P' \cap P'')=\rho(P') \cap \rho(P'')$, so it suffices to show that
$\sigma \subseteq \rho(P') \cap \rho(P'')$. Since
$\relint(\sigma) \cap \rho(P') \cap \rho(P'') \neq \emptyset$ and
$\rho(P')\in \Delta$, we have that $\sigma \subseteq \rho(P')$, so it
suffices to show $\sigma \subseteq \rho(P'')$.
Let $P \in \cP$ be a polyhedron containing $P'$. Then 
$\relint(\sigma) \cap \rho(P \cap P'') \neq \emptyset$, so by compatibility
$\sigma \subseteq \rho(P \cap P'') \subseteq \rho(P'')$, as desired.\qed

\begin{defn} Let $X$ and $X'$ be closed subschemes of $\T$, and fix a
  choice of polyhedral complex structures on $\Trop(X)$ and 
  $\Trop(X')$.  Let $C$ be a connected component of $\Trop(X)\cap\Trop(X')$. 
  A \emph{compactifying datum} for $X,X'$ and $C$ consists of a pair
  $(\Delta,\cP)$, where $\cP$ is  
  a finite collection of integral $G$-affine polyhedra in $N_\R$ such that 
  \[\Trop(X) \cap \Trop(X') \cap |\cP| = |C|, \]
  and $\Delta$ is an integral compactifying fan for $\cP$ which is compatible with
  $\Trop(X') \cap \cP$.
\end{defn}

The convention that $\Delta$ should be compatible
specifically with $\Trop(X') \cap \cP$ rather than either $\Trop(X) \cap \cP$
or $\Trop(X') \cap \cP$ is made out of convenience, to simplify the 
statements of Lemma~\ref{lem:moving.lemma} and 
Corollary~\ref{cor:moving-data} below.

\begin{rem} \label{rem:cP.is.C}
  If $\cP=C$ (with the induced polyhedral complex structure), then 
  in order for $(\Delta,\cP)$ to be a compactifying 
  datum for $X,X'$ and $C$, it suffices that $\Delta$ be an integral compactifying fan 
  for $C$, since such $\Delta$ is automatically compatible with $\Trop(X')\cap\cP$. In particular, by Remark~\ref{rem:cpct.fans.exist} compactifying data 
  always exist.  (The extra flexibility in the choice of $\cP$ will be used
  in the proof of Theorem~\ref{thm:main-multiple}.)
\end{rem}

\begin{lem}[Tropical moving lemma] \label{lem:moving.lemma}
  Let $X$ and $X'$ be closed subschemes of $\T$, and 
  suppose that $\codim(X) + \codim(X') = \dim(\T)$.  Choose polyhedral
  complex structures on $\Trop(X)$ and $\Trop(X')$.
  Let $C$ be a connected component of $\Trop(X)\cap\Trop(X')$ and
  let $(\Delta,\cP)$ be a compactifying datum for $X,X'$ and $C$.
  There exists a $\Delta$-thickening $\cP'$ of $\cP$, a number
  $\epsilon > 0$, and a cocharacter $v\in N$ with
  the following properties: 
  \begin{enumerate}
  \item $(\Delta,\cP')$ is a compactifying datum for $X$, $X'$ and $C$. 
  \item For all $r\in[-\epsilon,0)\cup(0,\epsilon]$, the set 
  $(\Trop(X) + r\cdot v)\cap\Trop(X') \cap |\cP'|$ is finite and contained in 
  $|\cP'|^{\circ}$, and each point lies in the interior of facets of
  $\Trop(X)+r\cdot v$ and $\Trop(X')$.
  \end{enumerate}
\end{lem}

\pf We begin with the observation that if $P,P'$ are disjoint
polyhedra then there exists a thickening of $P$ which is disjoint from
$P'$. Indeed, write 
$P = \bigcap_{i=1}^r \{ v\in N_\R~:~\angles{u_i,v}\leq a_i \}$
for $u_1,\ldots,u_r\in M$ and $a_1,\ldots,a_r\in G$, and for 
$t \geq 0$ set
$P_t = \bigcap_{i=1}^r \{ v\in N_\R~:~\angles{u_i,v}\leq a_i+t \}$.
Then $t \mapsto P_t\cap P'$ is a continuous family of
polyhedra with $P_0\cap P' = \emptyset$, so by
Corollary~\ref{cor:polyhedral-family} we have $P_t\cap P' = \emptyset$ for
some $t > 0$.

By Lemma~\ref{lem:delta.polyhedra} there exists a 
$\Delta$-decomposition $\cP''$ of $\cP$, and $\Delta$ is still compatible 
with $\cP'' \cap \Trop(X')$. Now, $\Delta$ is a compactifying fan for any
thickening $\cP'$ of $\cP''$. It follows from the 
above observation that $\cP'$ may be chosen such that
$\Trop(X) \cap \Trop(X') \cap |\cP'|=C$, and such that for each polyhedron
$P' \in \cP'$, if $P'$ is a thickening of $P'' \in \cP''$, then $P'$ meets 
precisely the same polyhedra of $\Trop(X')$ as $P''$. Given $P \in \Trop(X')$
meeting $P''$, note that $\rho(P \cap P'')=\rho(P) \cap \rho(P'') =
\rho(P) \cap \rho(P') = \rho(P \cap P')$, so the compatibility of $\Delta$
with $\Trop(X') \cap \cP'$ follows from the compatibility with 
$\Trop(X') \cap \cP''$. This proves~(1).

For any $v\in N$, in order to prove that there exists $\epsilon>0$ with
$(\Trop(X) + r\cdot v)\cap\Trop(X') \cap |\cP'| \subseteq |\cP'|^{\circ}$ for 
all $r\in[-\epsilon,0)\cup(0,\epsilon]$, we argue similarly to the above.
Indeed, note that
$r\mapsto (P + r\cdot v) \cap P'$ is a continuous family of polyhedra for
any polyhedra $P\subseteq\Trop(X),P'\subseteq\Trop(X')$, and that
$|\cP'|\setminus|\cP'|^\circ$ is contained in  
$|\cP'|\setminus\bigcup_{P\in\cP'} P^{\circ}$, which is a finite union of
polyhedra disjoint from $\Trop(X)\cap\Trop(X')$. The
finiteness assertion for suitable choice of $v$ follows from the fact that
$\dim(P)+\dim(P')\leq\dim(\T)$ for any polyhedra $P\subseteq\Trop(X)$ and
$P'\subseteq\Trop(X')$, since generic translates of
any two affine spaces of complementary dimension intersect in one or zero 
points. Similarly, any point lies in the interior of facets because the
lower-dimensional faces have subcomplementary dimension, and thus 
generic translates do not intersect.\qed

A tuple $(\cP', \epsilon, v)$
satisfying Lemma~\ref{lem:moving.lemma} will be called
a set of \emph{tropical moving data} for $(\Delta,\cP)$. 

\begin{cor}\label{cor:moving-data} In the situation of Lemma 
\ref{lem:moving.lemma}, we have 
\[ \bar{\Trop(X)} \cap \bar{\Trop(X')} \cap \bar{|\cP'|} = \bar{C}
\subseteq \bar{|\cP|} \subseteq \bar{|\cP'|}{}^{\circ},\]
and for all $r\in[-\epsilon,0)\cup(0,\epsilon]$ we have
\[\bar{\Trop(X) + r\cdot v}\cap\bar{\Trop(X')} \cap \bar{|\cP'|}= 
(\Trop(X) + r\cdot v)\cap\Trop(X') \cap|\cP'| \subseteq |\cP'|^{\circ}
\subseteq \bar{|\cP'|}{}^{\circ},\]
all closures being taken in $N_\R(\Delta)$.
\end{cor}

\pf This follows immediately from the compatibility hypotheses
of a compactifying datum, together with 
Proposition~\ref{prop:compactifying.fans}, noting that 
$\bar{|\cP|} \subseteq \bar{|\cP'|}{}^{\circ}$ by Remark 
\ref{rem:thickening.interior}, and, for the second statement, that the
closure of a finite set is itself.\qed

\paragraph[Relative boundary and properness in analytic geometry]
Our next goal is to construct a proper family of analytic spaces from
compactifying and moving data as above.  First we briefly review the
analytic notion of properness.
In this section we assume that $K$ is complete and nontrivially valued.
Recall from~\parref{par:analytic.spaces} that by an analytic space we mean a
Hausdorff, good, strictly $K$-analytic space.

Let $\cX\to\cY$ be a morphism of analytic spaces.  There is a canonical
open subset $\Int(\cX/\cY)$ of $\cX$ called the \emph{relative interior}
of the morphism $\cX\to\cY$ (not to be confused with the relative interior of a
polyhedron); its complement $\del(\cX/\cY)$ in $\cX$ is the  
\emph{relative boundary}.  The \emph{absolute interior} $\Int(\cX)$ of an
analytic space $\cX$ is the relative interior of the structure morphism
$\cX\to\sM(K)$, and the \emph{absolute boundary} $\del(\cX)$ is its
complement in $\cX$.  We will use
the following properties of the relative interior and relative boundary;
for the definition of $\Int(\cX/\cY)$
see~\cite[\S3.1]{berkovich:analytic_geometry}. 

\begin{prop}[Berkovich] \label{prop:relint.properties}
  \begin{enumerate}
  \item If $\cX$ is an analytic domain in an analytic space $\cY$ then
    $\Int(\cX/\cY)$ coincides with the topological interior of $\cX$ in $\cY$.
  \item Let $\cX\overset{f}\to\cY\to\cZ$ be a sequence of morphisms of analytic
    spaces.  Then 
    \[ \Int(\cX/\cZ) = \Int(\cX/\cY)\cap f\inv(\Int(\cY/\cZ)). \]
  \item Let $\cX\to\cY$ and $\cY'\to\cY$ be morphisms of analytic spaces,
    let $\cX' = \cY'\times_\cY\cX$, and let $f:\cX'\to\cX$ be the
    projection.  Then $f\inv(\Int(\cX/\cY))\subseteq\Int(\cX'/\cY')$.
  \item If $X$ is a finite-type $K$-scheme then 
    $\del(X^\an) = \emptyset$.
  \end{enumerate}
\end{prop}

See~\cite[Proposition~3.1.3 and Theorem~3.4.1]{berkovich:analytic_geometry} 
for the proofs.  The notion of a proper morphism of analytic spaces is
defined in terms of the relative interior:

\begin{defn}
  Let $f: \cX\to\cY$ be a separated morphism of analytic spaces.
  \begin{enumerate}
  \item $f$ is \emph{boundaryless} provided that
    $\del(\cX/\cY)=\emptyset$.
  \item $f$ is \emph{compact} provided that the inverse image of a compact
    set is compact.
  \item $f$ is \emph{proper} if it is both boundaryless and compact.
  \end{enumerate}
\end{defn}

Proper morphisms of analytic spaces behave much like proper morphisms of
schemes.  A morphism $X\to Y$ of finite-type $K$-schemes is proper if and only
if $X^\an\to Y^\an$ is proper.
A finite morphism $f: \cX\to\cY$ of analytic
spaces is proper.  (To say that $f$ is finite means that for every
affinoid domain $\sM(B)\subseteq\cY$, its inverse image 
$f\inv(\sM(B))$ is an affinoid domain $\sM(A)$, and $A$ is a finite $B$-module.)
The converse holds in the following familiar cases:

\begin{thm} \label{thm:finite.implies.proper}
  Let $f: \cX\to\cY$ be a proper morphism of analytic spaces.
  \begin{enumerate}
  \item If $\cX$ and $\cY$ are affinoid spaces then $f$ is finite.
  \item If $f$ has finite fibers then $f$ is finite.
  \end{enumerate}
\end{thm}

\pf Part~(1) follows from the Kiehl's direct image
theorem~\cite[Proposition~3.3.5]{berkovich:analytic_geometry}, and~(2) is
Corollary~3.3.8 of loc.\ cit.\qed

\paragraph[The tropical criterion for properness]
Fix an integral pointed fan $\Delta$ in $N_\R$.
Let $P$ be an integral $G$-affine polyhedron in $N_\R$ with 
$\rho(P)\in\Delta$, and let $\bar P$ be its closure in $N_\R(\Delta)$.  The
inverse image of $\bar P$ under $\trop:X(\Delta)^\an\to N_\R(\Delta)$ is
called a \emph{polyhedral domain} and is denoted $\cU_P$;
see~\cite[\S6]{jdr:trop_ps}.  This is an 
affinoid domain in $X(\Delta)^\an$.  If $\cP$ is a finite collection of
integral $G$-affine polyhedra with recession cones contained in $\Delta$,
then $\cU_\cP \coloneq \bigcup_{P\in\cP} \cU_P = \trop\inv(\bar{|\cP|})$
is a compact analytic domain in $X(\Delta)^\an$.  

\begin{lem} \label{lem:relint}
  Let $\Delta$ be an integral pointed fan in $N_\R$ and let $\cP$ be a
  finite collection of integral $G$-affine polyhedra with recession cones
  contained in $\Delta$.  Let $\cS$ be an analytic space and let 
  $p_2:\cS\times\cU_\cP\to\cU_\cP$ be the projection onto the second
  factor.  Then
  \[ \Int(\cS\times\cU_\cP/\cS) \supset 
  (\trop\circ p_2)\inv(\bar{|\cP|}{}^\circ), \]
  where the closure is taken in $N_\R(\Delta)$.
  In particular, $\Int(\cU_\cP)\supset\trop\inv(\bar{|\cP|}{}^\circ)$.
\end{lem}

\pf By Proposition~\ref{prop:relint.properties}(3) we have
$\Int(\cS\times\cU_\cP/\cS)\supset p_2\inv(\Int(\cU_\cP))$, so it suffices
to show that  $\Int(\cU_\cP)\supset\trop\inv(\bar{|\cP|}{}^\circ)$.
By Proposition~\ref{prop:relint.properties}(1), 
$\Int(\cU_\cP/X(\Delta)^\an)$ is the topological interior of $\cU_\cP$ in
$X(\Delta)^\an$ since $\cU_\cP$ is an analytic domain in $X(\Delta)^\an$.
Since $\trop:X(\Delta)^\an\to N_\R$ is continuous,
the set $\trop\inv(\bar{|\cP|}{}^\circ)$ is open in $X(\Delta)^\an$, so 
$\trop\inv(\bar{|\cP|}{}^\circ)\subseteq\Int(\cU_\cP/X(\Delta)^\an)$.
Applying Proposition~\ref{prop:relint.properties}(2) to the sequence of
morphisms $\cU_\cP\inject X(\Delta)^\an\to\sM(K)$, one obtains
\[ \Int(\cU_\cP) = \Int(\cU_\cP/X(\Delta)^\an) \cap \Int(X(\Delta)^\an). \]
But $\Int(X(\Delta)^\an) = X(\Delta)^\an$ by
Proposition~\ref{prop:relint.properties}(4), so
$\Int(\cU_\cP) = \Int(\cU_\cP/X(\Delta)^\an)\supset\trop\inv(\bar{|\cP|}{}^\circ)$.
\qed

\begin{lem} \label{lem:check.classical}
  Let $\Delta$ be an integral pointed fan in $N_\R$,
  let $\cP$ be a finite collection of integral $G$-affine polyhedra with
  recession cones contained in $\Delta$,
  let $\cS$ be an analytic space, and let 
  $\cX\subseteq\cS\times X(\Delta)^\an$ be a Zariski-closed subspace.
  Suppose that $\trop(\cX_s)\subseteq\bar{|\cP|}$ for all $s\in|\cS|$.  Then  
  $\cX\subseteq\cS\times\cU_\cP$.
\end{lem}

\pf The hypothesis in the statement of the lemma is equivalent to requiring
that $\cX_s\subseteq\{s\}\times\cU_\cP$ for all \emph{rigid} points
$s\in|\cS|$.  Since $|\cX|$ maps to $|\cS|$,
the set $\Djunion_{s\in|\cS|} \cX_s\supset|\cX|$ is everywhere
dense in $\cX$, so since $\Djunion_{s\in|\cS|} \cX_s$ is
contained in the closed subset $\cS\times\cU_\cP$, we have
$\cX\subseteq\cS\times\cU_\cP$.\qed

The following proposition can be found in~\cite[\S9]{jdr:trop_ps}, in a
weaker form and in the language of classical rigid spaces.

\begin{prop}[Tropical criterion for properness] \label{prop:properness}
  Let $\Delta$ be an integral pointed fan in $N_\R$,
  let $\cS$ be an analytic space,
  and let $\cX$ be a Zariski-closed subspace of 
  $\cS\times X(\Delta)^\an$.  Suppose that there exists a 
  finite collection $\cP$ of integral $G$-affine polyhedra with
  recession cones contained in $\Delta$ such that
  $\trop(\cX_s,\Delta)\subseteq\bar{|\cP|}$ for all $s\in|\cS|$,
  where the closure is taken in $N_\R(\Delta)$.
  Then $\cX\to\cS$ is proper.
  Moreover, if $\cP = \{P\}$ is a single polyhedron then $\cX\to\cS$ is finite.
\end{prop}

\pf By Lemma~\ref{lem:check.classical}, the condition on the
tropicalizations implies that $\cX\subseteq\cS\times\cU_\cP$,
i.e.\ that $\trop(p_2(\cX))\subseteq\bar{|\cP|}$.
Since properness can be checked affinoid-locally on the base, we may
assume that $\cS$ is affinoid.  Then $\cS\times\cU_\cP$ is compact, being
a finite union of affinoids, so $\cX$ is compact, and therefore
$\cX\to\cS$ is a compact map of topological spaces.
Replacing $\cP$ by a thickening, we can assume that 
$\trop(p_2(\cX))\subseteq\bar{|\cP|}{}^\circ$ (cf.\
Remark~\ref{rem:thickening.interior}), so 
$\cX\subseteq\Int(\cS\times\cU_\cP/\cS)$ by Lemma~\ref{lem:relint}.
Applying Proposition~\ref{prop:relint.properties}(2) to the sequence of
morphisms $\cX\inject\cS\times\cU_\cP\to\cS$ we obtain 
\[ \Int(\cX/\cS) = \Int(\cX/\cS\times\cU_\cP) \cap 
\Int(\cS\times\cU_\cP/\cS) = \Int(\cX/\cS\times\cU_\cP). \]
Since $\cX\to\cS\times\cU_\cP$ is a closed immersion it is finite, hence
proper, so $\cX = \Int(\cX/\cS\times\cU_\cP)$; therefore $\Int(\cX/\cS)=\cX$,
so $\cX\to\cS$ is boundaryless and compact, hence proper.

Now suppose that $\cP = \{P\}$ is a single polyhedron, still assuming $\cS$
affinoid.  Then $\cS\times\cU_\cP = \cS\times\cU_P$ is affinoid, so $\cX$
is affinoid; hence 
$\cX\to\cS$ is finite by Theorem~\ref{thm:finite.implies.proper}(1), being
a proper morphism of affinoids.\qed 

Since properness and finiteness can be checked after analytification, we
have the following algebraic consequence.  

\begin{cor}
  Let $\Delta$ be an integral pointed fan in $N_\R$ and let 
  $X\subseteq X(\Delta)$ be a closed subscheme.
  Suppose that there exists a 
  finite collection $\cP$ of integral $G$-affine polyhedra with
  recession cones contained in $\Delta$ such that
  $\Trop(X,\Delta)\subseteq\bar{|\cP|}$,
  where the closure is taken in $N_\R(\Delta)$.
  Then $X$ is proper.
  Moreover, if $\cP = \{P\}$ is a single polyhedron then $X$ is finite.
\end{cor}

\paragraph[The moving construction] \label{par:moving.construction}
Finally we show how a set of tropical moving data gives rise to a 
proper family over an analytic annulus.
Fix $X,X'\subseteq\T$ with $\codim(X)+\codim(X')=\dim(\T)$ and fix a
connected component $C$ of $\Trop(X)\cap\Trop(X')$.  Choose polyhedral
structures on $\Trop(X)$ and $\Trop(X')$,
let $(\Delta,\cP)$ be a compactifying datum for $X,X'$ and $C$, and
choose a set $(\cP',\epsilon,v)$ of tropical moving data for $(\Delta,\cP)$.
We may assume without loss of generality that $\epsilon\in G$.  Let
\[ \cS_\epsilon = \cU_{[\epsilon,\epsilon]} =
\val\inv([-\epsilon,\epsilon]) \subset\G_m^\an; \]
this is the annulus whose set of
$\bar K$-points is $\{t\in \bar K{}^\times~:~\val(t)\in[-\epsilon,\epsilon]\}$.  
It is a polytopal domain (and in particular an affinoid domain) in $\G_m^\an$.  

Let $\bar X$ and $\bar X\p$ denote the closures of $X$ and $X'$ in $X(\Delta)$,
respectively.  Considering $v$ as a homomorphism $v:\G_m\to\T$,
we obtain an action $\mu:\G_m\times X(\Delta)\to X(\Delta)$ 
of $\G_m$ on $X(\Delta)$ given by $\mu(t,x) = v(t)\cdot x$.  Note that 
$(p_1,\mu):\G_m\times X(\Delta)\to\G_m\times X(\Delta)$ is an
isomorphism, where $p_1$ is projection onto the first factor.  Let 
$\fX\coloneq (p_1,\mu)(\G_m\times\bar X)$ and
$\fX'\coloneq \G_m\times\bar X\p$.
These are closed subschemes of $\G_m\times X(\Delta)$, which we will think
of as being flat families of closed subschemes of $X(\Delta)$ parameterized by
$\G_m$.  A point $t\in\G_m^\an$ can be thought of as a morphism
$t: \sM(\sH(t))\to(\G_m^\an)_{\sH(t)}$, which is given by an element of
$\sH(t)^\times$, and is thus the analytification of a
morphism $t: \Spec(\sH(t))\to(\G_m)_{\sH(t)}$.  Since analytifications commute with
fiber products and extension of scalars, the fiber $\fX_t$ of $\fX^\an$ over $t$ is
naturally identified with the analytification of 
$v(t)\cdot\bar X_{\sH(t)}$, which is the closure of 
$v(t)\cdot X_{\sH(t)}$ in $X(\Delta)_{\sH(t)}$.

Let $\fY = \fX\cap\fX'\subseteq\G_m\times X(\Delta)$ and let 
\[ \cY = \fY^\an\cap(\cS_\epsilon\times\cU_{\cP'})
= \fY^\an\cap(\val\circ p_1)\inv([-\epsilon,\epsilon])\cap
(\trop\circ p_2)\inv(\bar{|\cP'|}). \]
This is a Zariski-closed subspace of
$\cS_\epsilon\times\cU_{\cP'}$.
For $t\in\cS_\epsilon$ we have
\[ \Trop\big(v(t)\cdot\bar X_{\sH(t)},\Delta\big) = 
\Trop(\bar X,\Delta) - \val(t)\cdot v. \]

\begin{prop} \label{prop:conn.comp}
  The analytic space $\cY$ is a
  union of connected components of  
  $\fY^\an \cap (\cS_\epsilon\times X(\Delta)^\an)$.  Moreover, $\cY$ is
  proper over $\cS_\epsilon$ and Zariski-closed in 
  $\cS_\epsilon\times X(\Delta)^\an$.  
\end{prop}

\pf Since $\cY$ is the intersection of 
$\fY^\an\cap (\cS_\epsilon\times X(\Delta)^\an)$ with the (compact)
affinoid domain $\cS_\epsilon\times\cU_{{\cP'}}$, it is closed.
On the other hand, it follows from Corollary~\ref{cor:moving-data} and
Lemma~\ref{lem:check.classical} that 
$\trop(p_2(\cY))\subseteq\bar{|\cP'|}{}^\circ$, where 
$p_2: \cS_\epsilon\times X(\Delta)^\an\to X(\Delta)^\an$ is projection onto 
the second factor. Thus, $\cY$ is the intersection of 
$\fY^\an\cap (\cS_\epsilon\times X(\Delta)^\an)$ with the open subset
$(\trop\circ p_2)\inv(\bar{|{\cP'}|}{}^\circ)$, so $\cY$ is both open and
closed in $\fY^\an\cap (\cS_\epsilon\times X(\Delta)^\an)$.
Hence $\cY$ is Zariski-closed in $\cS_\epsilon\times X(\Delta)^\an$, so
$\cY\to\cS_\epsilon$ is proper by Proposition~\ref{prop:properness}.\qed

\section{Continuity of intersection numbers}\label{sec:int.nums}

In this section we prove a ``continuity of intersection numbers'' theorem
in the context of a relative dimension-zero intersection of flat families
over an analytic base.  We will apply this in section~\ref{sec:mainthm} to
the family constructed in~\parref{par:moving.construction}.

In this section we assume that $K$ is complete and nontrivially valued.

\paragraph[Flat and smooth morphisms of analytic spaces]
We begin with a review of flatness and smoothness in
analytic geometry.
In general the notion of a flat morphism of analytic spaces is quite
subtle.  However, since we are assuming that all of
our analytic spaces are strictly $K$-analytic, separated, and good, the
situation is much simpler: a morphism $f:\cY\to\cX$ of analytic spaces is
\emph{flat} provided that, for every pair of affinoid domains
$\cV = \sM( B)\subseteq\cY$ and $\cU=\sM(A)\subseteq\cX$ with
$f(\cV)\subseteq\cU$, the corresponding homomorphism $A\to B$ is flat;
see~\cite[Corollary~7.2]{ducros:flatness} (Ducros calls this notion
``universal flatness'').  This condition can be
checked on an affinoid cover.  

The notion of smoothness  that is relevant
for our purposes is called ``quasi-smoothness'' by Ducros in loc.\ cit.\
and ``rig-smoothness'' in the language of classical rigid spaces.  A
morphism $f:\cY\to\cX$ is said to be \emph{quasi-smooth} if it is flat
with geometrically regular
fibers~\cite[Proposition~3.14]{ducros:flatness}.  

Both flatness and quasi-smoothness are
preserved under composition and change of base, and the inclusion of
an analytic domain is flat and quasi-smooth.  
A morphism $Y\to X$ of finite-type $K$-schemes is flat
(resp.\ smooth) if and only if $Y^\an\to X^\an$ is flat (resp.\
quasi-smooth).  

\begin{rem}
  The best reference for the notions of flatness and smoothness in
  Berkovich's language is~\cite{ducros:flatness}; however, 
  Ducros works in much greater generality than is necessary for our
  purposes.  Most of the results that we will use have been
  known for much longer, but can only be found in the literature in the
  language of classical rigid spaces.
\end{rem}

We define local intersection numbers of schemes and
analytic spaces using a modification of Serre's definition:

\begin{defn} \label{defn:int.num}
  Let $Y$ be a smooth scheme over a field $k$ (resp.\ a quasi-smooth
  analytic space over a nontrivially valued complete non-Archimedean field
  $k$), let $X,X'\subseteq Y$ 
  be closed subschemes (resp.\ Zariski-closed subspaces), and suppose that
  $x\in|X\cap X'|$ is an   isolated point of $X\cap X'$.  The 
  \emph{local intersection number} of $X$ and $X'$ at $x$ is defined to be
  \[ i_k(x,\, X\cdot X';\, Y) = \sum_{i=0}^{\dim(Y)} 
  (-1)^i \dim_k \Tor_i^{\sO_{Y,x}} (\sO_{X,x},\,\sO_{X',x}). \]
  If $X\cap X'$ is $K$-finite,
  the \emph{intersection number} of $X$ and $X'$ is 
  \[ i_k(X\cdot X';\, Y) 
  = \sum_{x\in|X\cap X'|} i_k(x,\, X\cdot X';\, Y). \]
\end{defn}

\begin{rem} \label{rem:vanish.tor}
  The dimension of $X\cap X'$ is zero at an isolated point $x$ of
  $X\cap X'$.  Hence $\sO_{X\cap X',x}$ 
  is an Artin local ring, being Noetherian of Krull dimension zero.  The
  finitely generated $\sO_{Y,x}$-module 
  $\Tor_i^{\sO_{Y,x}} (\sO_{X,x},\,\sO_{X',x})$ is naturally an
  $\sO_{X\cap X',x}$-module, and is therefore finite-dimensional
  over $k$.  Moreover, $\sO_{Y,x}$ is a regular local ring 
  as $Y$ is smooth (resp.\ quasi-smooth) over $k$; hence we have
  $\Tor_i^{\sO_{Y,x}}(\sO_{X,x},\,\sO_{X',x}) = 0$ for $i > \dim(Y)$.

\end{rem}

\begin{rem} \label{rem:global.int.nums}
  Suppose that $X\cap X'$ is finite.  The coherent sheaf 
  $\sTor_i^{\sO_Y}(\sO_X,\sO_{X'})$ is supported on $X\cap X'$.  Hence its
  space of global sections 
  $\Tor_i^{\sO_Y}(\sO_X,\sO_{X'}) = \Gamma(Y,\sTor_i^{\sO_Y}(\sO_X,\sO_{X'}))$
  breaks up as
  \[ \Tor_i^{\sO_Y}(\sO_X,\sO_{X'}) 
  = \Dsum_{x\in|X\cap X'|} \Tor_i^{\sO_{Y,x}} (\sO_{X,x},\,\sO_{X',x}), \]
  so it follows that 
  \[ i_k(X\cdot X';\, Y) = 
  \sum_{i=0}^{\dim(Y)} (-1)^i \dim_k \Tor_i^{\sO_Y}(\sO_X,\sO_{X'}). \]
  Hence our definition agrees
  with~\cite[Definition~4.4.1]{osserman_payne:lifting}. 
\end{rem}

\begin{rem}
  It is clear that $i_k(x,X\cdot X';Y)$ is local on $Y$, in that it only
  depends on an affine (resp.\ affinoid) neighborhood of $x$.
\end{rem}

We have the following compatibility of algebraic and analytic intersection
numbers: 

\begin{prop} \label{prop:alg.an.int.nums}
  Let $Y$ be a smooth scheme over $K$, let $X,X'\subseteq Y$ be closed
  subschemes, and let $x\in |X\cap X'|$ be an isolated point of 
  $X\cap X'$.  Then
  \[ i_K(x,\, X\cdot X';\, Y) = i_K(x,\, X^\an\cdot(X')^\an;\, Y^\an) \]
  under the identification of $|X\cap X'|$ with $|X^\an\cap(X')^\an|$.
\end{prop}

\pf By~\cite[Lemma~A.1.2(2)]{conrad:irredcomps} the local rings
$\sO_{Y,x}$ and $\sO_{Y^\an,x}$ have the same completion, so the proposition
follows from Lemma~\ref{lem:extend.tor} below.\qed

Our goal will be to prove the following invariance of  intersection
numbers in families over analytic spaces:

\begin{prop} \label{prop:cont_nums}
  Let $\cS$ be an analytic space, let 
  $\cZ$ be a quasi-smooth analytic space, and let
  $f:\cZ\to\cS$ be a quasi-smooth morphism.
  Let $\cX,\cX'\subseteq\cZ$ be Zariski-closed subspaces, flat over $\cS$,
  such that $\cY = \cX\cap\cX'$ is finite over $\cS$.  Then
  the map
  \[ s~\longmapsto~ i_{K(s)}(\cX_s\cdot\cX'_s;\, \cZ_s)\quad:\quad
  |\cS|\To\Z \]
  is constant on connected components of $\cS$.
\end{prop}

\smallskip 

We will need the following lemmas.

\begin{lem} \label{lem:extend.tor}
  Let $A$ be a noetherian local ring with maximal ideal $\fm$ and let 
  $\hat A$ be its $\fm$-adic completion.  Let $M,N$ be
  finitely generated $A$-modules such that 
  $\Supp(M)\cap\Supp(N) = \{\fm\}$.  Then for all $i\geq 0$,
  the natural map
  \[ \Tor_i^A(M,\,N) \To 
  \Tor_i^{\hat A}(M\tensor_A\hat A,\, N\tensor_A\hat A) \]
  is an isomorphism.
\end{lem}

\pf Let $\fa = \Ann(M) + \Ann(N)$.
Since $\fm/\fa$ is a nilpotent ideal in $A/\fa$, the finitely
generated $A/\fa$-module $\Tor_i^A(M,N)$ is $\fm$-adically discrete, so 
$\Tor_i^A(M,N)\to\Tor_i^A(M,N)\tensor_A\hat A$ is an isomorphism.  But
$A\to\hat A$ is flat, so 
\[ \Tor_i^{\hat A}(M\tensor_A\hat A,\, N\tensor_A\hat A) \cong
\Tor_i^A(M,N)\tensor_A\hat A \]
naturally.\qed

Recall from~\parref{par:analytic.spaces} that by an analytic space we mean a
Hausdorff, good, strictly $K$-analytic space.

\begin{lem} \label{lem:finite.in.affinoid}
  Let $f:\cZ\to\cS$ be a morphism of analytic spaces and let
  $\cY\subseteq\cZ$ be a Zariski-closed subspace which is finite over
  $\cS$.  Then for any point $s\in\cS$, there exists an affinoid neighborhood
  $\cU$ of $s$ and an affinoid domain $\cV\subseteq f\inv(\cU)$ such that
  $\cY\cap f\inv(\cU)\subseteq\cV$.
\end{lem}

\pf Fix $s\in\cS$.  We may replace $\cS$ with an affinoid neighborhood of
$s$ to assume $\cS$ affinoid.  For $y\in\cY_s\coloneq f\inv(s)\cap\cY$ 
let $\cV(y)$ be an
affinoid neighborhood of $y$ in $\cZ$.  We may choose the $\cV(y)$ such that
$\cV(y)\cap\cV(y') = \emptyset$ for $y\neq y'$; this is possible because
$\cY_s$ is a finite set of points in the Hausdorff space $\cZ$,
and the affinoid neighborhoods of a point form a base of closed
neighborhoods around that point.  Let $\cV(y)^\circ$ denote the interior
of $\cV(y)$ in $\cZ$ and let
$C = f(\cY\setminus\bigcup_{y\in \cY_s} \cV(y)^\circ)$.
Since $\cS$ is affinoid and $\cY\to\cS$ is finite, 
$\cY$ is affinoid, hence compact; therefore 
$\cY\setminus\bigcup_{y\in \cY_s} \cV(y)^\circ$ is compact,
so $C$ is compact, hence closed in $\cS$.  By construction, a point
$s'\in\cS$ is not contained in $C$ if and only if 
$\cY_{s'}\subseteq\bigcup_{y\in\cY_s}\cV(y)^\circ$; in particular, 
$s\notin C$.  Let $\cU$ be an affinoid neighborhood of $s$ contained in
$\cS\setminus C$ and let 
$\cV = f\inv(\cU)\cap\bigcup_{y\in\cY_s} \cV(y)$.  Clearly
$f\inv(\cU)\cap\cY\subseteq\cV$.  Since the $\cV(y)$ are disjoint, 
the union $\bigcup_{y\in\cY_s} \cV(y)$ is affinoid, so $\cV$ is affinoid,
being a fiber product of affinoids.\qed

\pf[of Proposition~\ref{prop:cont_nums}]
The question is local on $\cS$, in the following sense.  The analytic
space $\cS$ is connected if and only if the associated classical
rigid-analytic space $|\cS|$ is connected --- in other words, if and only
if the set $|\cS|$ is connected with respect to the Grothendieck topology
generated by subsets of the form $|\cU|$ for $\cU\subseteq\cS$ affinoid,
with coverings being the so-called admissible coverings.%
\footnote{It would be more elegant to prove the proposition for all points
  $s\in\cS$, but in that case one would have to treat the issue of
  exactness of the completed tensor product.}
Concretely, this means that if we can cover $\cS$ by 
affinoid domains $\{\cS_i\}_{i\in I}$ such that every point of $\cS$ is
contained in the interior of some $\cS_i$, then it suffices to prove the
proposition after base change to each $\cS_i$.
See~\cite[\S1.6]{berkovich:etalecohomology}. 
By Lemma~\ref{lem:finite.in.affinoid} we may assume that $\cS = \sM(R)$ and
$\cZ = \sM(C)$ are affinoid, and that $\cS$ is connected.  
Hence $\cX = \sM(A)$, $\cX' = \sM(A')$, and $\cY = \sM(B)$ affinoid as
well. 

Now we proceed as in the proof
of~\cite[Theorem~4.4.2]{osserman_payne:lifting}. 
Let $P_\bullet$ be a resolution of $A$ by finite free $C$-modules and let
$Q_\bullet = P_\bullet\tensor_C A'$.  Then the homology of $Q_\bullet$
calculates the groups $\Tor^C_i(A,A')$.  
For any maximal ideal $\fp\subset C$ the localization of 
$\Tor^C_i(A,A')$ at $\fp$ is canonically isomorphic to
$\Tor^{C_\fp}_i(A_\fp,A'_\fp)$; since $\cZ$ is quasi-smooth, 
its local rings are regular, so $C_\fp$ is
regular~\cite[Proposition~7.3.2/8]{bgr:nonarch}, and hence 
$\Tor^{C_\fp}_i(A_\fp,A'_\fp) = 0$  for $i > \dim(\cZ)$.  It follows that
$\Tor^C_i(A,A') = 0$ for $i > \dim(\cZ)$.

Let $s\in|\cS|$ and write $B_s = B\tensor_R K(s)$,
$C_s = C\tensor_R K(s)$, etc.
Since $A$ is $R$-flat,
$P_\bullet\tensor_R K(s)$ is a resolution of $A_s$
by finite free $C_s$-modules, so
$Q_\bullet\tensor_R K(s)$ computes
$\Tor_\bullet^{C_s}(A_s,A'_s)$.
Let $y$ be a point of $\cY_s$, and let $\fp$ be the corresponding maximal
ideal of $C_s$.  By~\cite[Proposition~7.3.2/3]{bgr:nonarch}, the local rings
$(C_s)_\fp$ and $\sO_{\cY_s,y}$ have the same completion.
Since $\Tor_i^{C_s}(A_s,A'_s)$
is supported on the finite set of points of $\cY_s$ for all $i$, we have
\[ \Tor_i^{C_s}(A_s,A'_s)
\cong \Dsum_{\fp\in\cY_s}
\Tor_i^{(C_s)_\fp} \big((A_s)_\fp, (A'_s)_\fp\big)
\cong \Dsum_{y\in|\cY_s|} 
\Tor_i^{\sO_{\cY_s,y}} \big(\sO_{\cX_s,y}, \sO_{\cX'_s,y}\big), \]
where the last equality comes from Lemma~\ref{lem:extend.tor}.
Therefore
\[ \sum_{i=0}^\infty (-1)^i \dim_{K(s)}\Tor_i^{C_s}(A_s,A'_s)
= i_{K(s)}(\cX_s\cdot\cX'_s;\, \cZ_s). \]

The finite $C$-modules $\Tor^C_\bullet(A,A')$ are supported on $\cY$, so
since $\cY$ is finite over $\cS$, they are in fact finite $R$-modules.
Viewing $Q_\bullet$ as a complex of $R$-modules with finitely many
$R$-finite cohomology groups, it follows 
from~\cite[Corollaire~0.11.9.2]{egaiii_1}
that there exists a quasi-isomorphic bounded below complex $M_\bullet$ of
free $R$-modules of finite (constant) rank.  Furthermore $Q_\bullet$ is a
complex of finite free $A'$-modules, hence flat $R$-modules, so by
Remark~11.9.3 of~loc.\ cit., for $s\in|\cS|$ the complex
$M_\bullet\tensor_R K(s)$ computes the homology of 
$Q_\bullet\tensor_R K(s)$, i.e.\ the groups
$\Tor_\bullet^{C_s}(A_s,A'_s)$.
Therefore
\[\begin{split} i_{K(s)}(\cX_s\cdot\cX'_s;\, \cZ_s)
&= \sum_i (-1)^i \dim_{K(s)}(\Tor_i^{C_s}(A_s,A'_s))\\
&= \sum_i (-1)^i \dim_{K(s)}(M_i\tensor_R K(s))
= \sum_i (-1)^i \rank_R(M_i) 
\end{split}\]
is independent of $s\in|\cS|$.\qed

\section{Tropical lifting theorems}\label{sec:mainthm}

We are now in a position to prove the main theorems relating algebraic and
tropical intersection multiplicities when the algebraic intersection is
finite but the tropical intersection is not necessarily finite (i.e.\ the
tropicalizations do not meet properly in the terminology
of~\cite{osserman_payne:lifting}).  First we prove the theorem for
intersections of two subschemes, then we extend to intersections of
several subschemes. In this section, we assume that $K$ is a possibly
trivially-valued non-Archimedean field which is complete or algebraically closed.

\paragraph[Tropical intersection multiplicities]
We begin by recalling the basic definitions of tropical intersection
theory. Let $X,X'\subseteq\T$ be pure-dimensional closed subschemes such that
$\codim(X)+\codim(X')= \dim(\T)$. We say that $\Trop(X)$ and $\Trop(X')$ intersect
\emph{tropically transversely} at a point $v \in \Trop(X) \cap \Trop(X')$
if $v$ is isolated and lies in the interior of facets in both $\Trop(X)$ and
$\Trop(X')$. If $\Trop(X)$ and $\Trop(X')$ intersect tropically transversely
at $v$, then the \emph{local tropical
intersection multiplicity} $i(v,\Trop(X) \cdot \Trop(X'))$ is defined
to be $[N:N_{ P}+N_{ P'}] m( P) m( P')$ where
$ P, P'$ are the facets of $\Trop(X)$ and $\Trop(X')$ respectively
containing $v$, we denote by $N_{ P}$ (respectively, $N_{ P'}$)
the sublattice of $N$ spanned by the translation of $ P$ (respectively,
$ P'$) to the origin, and $m( P)$ (respectively, $m( P')$)
denotes the multiplicity of $ P$ in $\Trop(X)$ (respectively,
of $ P'$ in $\Trop(X')$).

Now, suppose that $\Trop(X)$ does not meet $\Trop(X')$ tropically
transversely. The theory of stable tropical intersection allows us to
nonetheless define nonnegative intersection multiplicities at all points 
of $\Trop(X) \cap \Trop(X')$, such that the multiplicities will be positive
at only finitely many points. As in Lemma \ref{lem:moving.lemma},
for a fixed generic cocharacter $w \in N$, and 
sufficiently small $t>0$, we will have that $(\Trop(X)+tw)$ intersects
$\Trop(X')$ tropically transversely. Moreover, for $t$ sufficiently small,
which facets of $(\Trop(X)+tw)$ and $\Trop(X')$ meet one another is
independent of $t$. For $v \in \Trop(X) \cap \Trop(X')$, we can thus
define the \emph{local tropical intersection multiplicity} to be
\begin{equation}\label{eq:stable-trop-int}
i\big(v,\,\Trop(X) \cdot\Trop(X')\big) = \sum_{ P \ni v,  P' \ni v}
i\big(( P+tw)\cap P',\,(\Trop(X)+tw)\cdot \Trop(X')\big),\end{equation}
where as above $ P$ and $ P'$ are facets of $\Trop(X)$ and $\Trop(X')$
respectively. The fact that this definition is independent of the choice of 
$w$ is a consequence of the balancing condition for tropicalizations,
or can be seen 
algebraically via the close relationship to the intersection theory of
toric varieties; see \cite{fulton_sturmfels:itheory}.

For us, the relevant properties of \eqref{eq:stable-trop-int} are the
following, which are easy consequences of the definition:

\begin{prop}\label{prop:trop-int-basic} Let $X,X'$ be pure-dimensional
closed subschemes of $\T$ of complementary codimension.
\begin{enumerate} 
\item If $\Trop(X)$ intersects $\Trop(X')$ tropically transversely at $v$, 
then the two definitions above of $i(v,\Trop(X),\Trop(X'))$ agree.
\item In general, if $v \in \Trop(X) \cap \Trop(X')$, and $\epsilon>0$,
there exists $\delta >0$ such that for all $t<\delta$, every point
$( P+tw)\cap  P'$ occurring in \eqref{eq:stable-trop-int} is
within $\epsilon$ of $v$.
\end{enumerate}
\end{prop}

With these preliminaries out of the way, our starting point is the theorem 
of Osserman and Payne which guarantees
the compatibility of local tropical intersection multiplicities with
local algebraic intersection multiplicities, when the tropicalizations
intersect properly. The following is a special case
of~\cite[Theorem~5.1.1]{osserman_payne:lifting} (note that the
hypothesis that $X,X'$ are subvarieties is not used anywhere in the
proof).

\begin{thm} \label{thm:compat.proper}
  Suppose that $K$ is algebraically closed.
  Let $X,X'\subseteq\T$ be pure-dimensional closed subschemes of
  complementary codimension.  Let 
  $v\in\Trop(X)\cap\Trop(X')$ be an isolated point.  Then there are
  only finitely many points $x\in|X\cap X'|$ with $\trop(x)=v$, and
  \[ \sum_{\substack{x\in|X\cap X'|\\\trop(x)=v}} 
  i_K\big(x,\,X\cdot X';\,\T\big) = i\big(v,\,\Trop(X)\cdot\Trop(X')\big). \]
\end{thm}

We extend the above theorem to a higher-dimensional connected component of
$\Trop(X)\cap\Trop(X')$ as follows.

\begin{thm} \label{thm:mainthm}
  Let $X,X'$ be pure-dimensional closed subschemes of $\T$ of
  complementary codimension.
  Choose polyhedral complex structures on
  $\Trop(X)$ and $\Trop(X')$.  Let $C$ be a connected component of
  $\Trop(X)\cap\Trop(X')$ and let $(\Delta,\cP)$ be a compactifying datum 
  for $X,X'$ and $C$ such that
  $X(\Delta)$ is smooth.  Let $\bar X,\bar X\p$ be the closures of $X,X'$
  in $X(\Delta)$, respectively, and let $\bar C$ be the closure of $C$ in
  $N_\R(\Delta)$.  If there are only finitely many points 
  $x\in|\bar X\cap\bar X\p|$ with $\trop(x)\in\bar C$ then
  \begin{equation} \label{eq:mainthm}
    \sum_{\substack{x\in|\bar X\cap\bar X\p|\\\trop(x)\in\bar C}}
    i_K\big(x,\, \bar X\cdot\bar X\p;\, X(\Delta)\big) =
    \sum_{v\in C} i\big(v,\, \Trop(X)\cdot\Trop(X')\big). 
  \end{equation}
\end{thm}

\pf Let $K'$ be a complete, nontrivially valued, algebraically closed 
valued field extension of $K$. We claim that it is enough to prove the
theorem after extending to $K'$.
Since the weights on $\Trop(X)$ and $\Trop(X')$ are insensitive to
valued field extensions, the same is true for the tropical intersection
multiplicities, so we need only show that the appropriate sums of the local 
algebraic intersection multiplicities are preserved after extending scalars 
to $K'$. Given $v \in \bar C$, let $U$ be an open subscheme of 
$X(\Delta)$ such that $|\bar X\cap \bar X\p \cap U|$ is the set of all 
points of $|\bar X\cap \bar X\p|$ tropicalizing to $v$.  Then
\begin{equation} \label{eq:int.nums.tor}
  \sum_{\substack{x\in|\bar X\cap \bar X\p|\\\trop(x)=v}} 
i_K\big(x,\,\bar X\cdot \bar X\p;\,X(\Delta)\big) 
= i_K\big((\bar X\cap U)\cdot(\bar X\p \cap U);\, U\big)
= \sum_{i=0}^{\dim(\T)} (-1)^i \dim_K 
\Tor_i^{\sO_U} (\sO_{\bar X\cap U},\sO_{\bar X\p\cap U}). 
\end{equation}
The right side of the above equation is visibly insensitive to field
extensions.  Since 
$\bar X_{K'}\cap \bar X\p_{K'}\cap U_{K'} 
= (\bar X\cap \bar X\p\cap U)_{K'}$ is a finite
subscheme and $\trop(x) = v$ for all $x\in|(\bar X\cap \bar X\p\cap U)_{K'}|$,
we can apply~\eqref{eq:int.nums.tor} again after extending scalars to 
obtain the desired compatibility. We thus replace $K$ with $K'$, and
assume that $K$ is both algebraically closed and complete with respect
to a nontrivial valuation.

Now, let $({\cP'},\epsilon,v)$ be a set of tropical moving data for
$(\Delta,\cP)$ as in Lemma~\ref{lem:moving.lemma}, with $\epsilon\in G$.
Let $Z = \G_m\times X(\Delta)$, and let
$\fX,\fX'\subseteq Z$ be the closed subschemes defined
in~\parref{par:moving.construction}.  Carrying out the construction
of~\parref{par:moving.construction}, we let 
$\cS = \cS_\epsilon\subset\G_m^\an$, 
$\fY = \fX\cap\fX'$, and $\cY = \fY^\an\cap(\cS\times\cU_{\cP'})$.  By
Proposition~\ref{prop:conn.comp}, $\cY$ is a union of connected
components of $\fY^\an\cap(\cS\times X(\Delta)^\an)$, $\cY$ is
Zariski-closed in $\cS\times X(\Delta)^\an$, and $\cY\to\cS$ is
proper.  For $s\in\cS$ let
$\cY_s$ be the fiber of $\cY$ over $s$, and let
$T = \{ s\in\cS~:~\dim(\cY_s) > 0 \}$.
By construction the fiber of $\cY$ over $1\in|\G_m|$ is equal
$\{y\in(\bar X\cap\bar X\p)^\an~:~\trop(y)\in\bar C\}$, and by hypothesis
$1\notin T$.  The theorem on semicontinuity
of fiber dimension of morphisms of analytic
spaces~\cite[Theorem~4.9]{ducros:relative_dimension} then 
gives that $T$ is a finite set of rigid points of $\cS$.  Replacing $\cS$
with $\cS\setminus T$, we have that $\cY\to\cS$ is finite by
Theorem~\ref{thm:finite.implies.proper}(2). 

Applying Proposition~\ref{prop:cont_nums} with 
$\cZ = \cS\times\cU_{\cP'}$, $\cX = \fX^\an\cap\cZ$, and
$\cX' = (\fX')^\an\cap\cZ$, and using the compatibility of analytic and
algebraic local intersection numbers from
Proposition~\ref{prop:alg.an.int.nums}, we obtain that for all
$s\in|\cS|$,
\[\begin{split}
  \sum_{\substack{y\in|\bar X\cap\bar X\p|\\\trop(y)\in\bar C}} 
  &i_{K}\big(y,\, \bar X\cdot\bar X\p;\, X(\Delta)\big) \\
  &=  \sum_{\substack{y\in|(v(s)\cdot\bar X_{K(s)})\cap\bar X\p_{K(s)}|\\\trop(y)\in\bar{|\cP'|}}}
i_{K(s)}\big(y,\, (v(s)\cdot\bar X_{K(s)})\cdot\bar X\p_{K(s)};\, 
X(\Delta)_{K(s)}\big). 
\end{split}\]
But if $\val(s)\neq 0$ then $\Trop(X) - \val(s)\cdot v$ meets
$\Trop(X')$ properly, so by
Theorem~\ref{thm:compat.proper}, the right side of the
above equation is equal to 
\[ \sum_{v\in|\cP|} i\big(v,\,(\Trop(X)-\val(s)\cdot v)\cdot\Trop(X')\big).  \]
It follows from Proposition~\ref{prop:trop-int-basic} that for
$\val(s)$ sufficiently close to $0$, the above quantity is equal
to $\sum_{v\in C} i(v,\Trop(X)\cdot\Trop(X'))$,
which finishes the proof.\qed 

\begin{rem}
  Recall from Remark~\ref{rem:cP.is.C} that 
  if $\cP$ is the polyhedral
  complex underlying $C$ and $\Delta$ is any integral 
  compactifying fan for $C$, then $(\cP,\Delta)$ is a compactifying datum
  for $X,X',$ and $C$.  We will use the more general statement of
  Theorem~\ref{thm:mainthm} when proving Theorem~\ref{thm:main-multiple}
  below. 
\end{rem}

Proposition~\ref{prop:single.polyhedron} below is used to prove
Proposition~\ref{prop:automatic.finiteness}, which guarantees the
finiteness hypothesis in Theorems~\ref{thm:mainthm}
and~\ref{thm:main-multiple}.  See Remark~\ref{rem:fan.existence.3}.

\begin{prop} \label{prop:single.polyhedron}
  Let $\cP$ be a finite collection of integral $G$-affine 
  polyhedra in $N_\R$, and suppose 
  that $\Delta$ is an integral compactifying fan for $\cP$. Then 
  there exists an integral $G$-affine pointed polyhedron $P$ with 
  $|\cP| \subseteq P$ and $\rho(P) \in \Delta$ if and only if
  there exists $\sigma \in \Delta$ such that for all 
  $P' \in \cP$, the cone $\rho(P')$ is a face of $\sigma$. 
\end{prop}

\pf First suppose that $P$ exists as in the statement, 
and set $\sigma=\rho(P)$. We
then have to show that $\rho(P')\prec\sigma$ for all $P' \in \cP$.
Clearly $\rho(P') \subseteq \rho(P)=\sigma$, and $\rho(P')$ is a 
union of cones in $\Delta$, each of which must then be a face of $\sigma$.
We then conclude that $\rho(P')$ is a face of $\sigma$, as desired.
Conversely, suppose that we have $\sigma \in \Delta$ as in the statement.
Noting that $\cP$ consists entirely of pointed polyhedra, let $V$ be the
set of all vertices of all polyhedra $P'\in\cP$.  The convex
hull $\conv(V)$ of $V$ is a polytope, and any pointed polyhedron is the 
Minkowski sum of its recession cone and the convex hull of its vertices, so 
$|\cP|\subseteq P\coloneq \conv(V)+\sigma$.  This $P$ is an integral 
$G$-affine polyhedron with recession cone $\sigma$.\qed

\begin{prop} \label{prop:automatic.finiteness}
  In the situation of Theorem~\ref{thm:mainthm}, suppose in addition that
  the equivalent conditions of Proposition~\ref{prop:single.polyhedron}
  are satisfied for $\Delta$ and $\cP$.  Then
  there are automatically only finitely many points $x\in|\bar X\cap\bar X\p|$
  with $\trop(x)\in\bar C$.
\end{prop}

\pf As in the proof of Theorem \ref{thm:mainthm}, it is clearly enough to
consider the case that $K$ is complete, with nontrivial valuation.
Define $(\cP',\epsilon,v)$, and $\cY$ as in the proof of
Theorem~\ref{thm:mainthm}.  Then
\[ \cY_1 = (\bar X\cap\bar X\p)^\an\cap\cU_{\cP'} 
= (\bar X\cap\bar X\p)^\an\cap\cU_{\cP} 
= \{y\in(\bar X\cap\bar X\p)^\an~:~\trop(y)\in\bar C\} \]
is Zariski-closed in $X(\Delta)^\an$.  By
Proposition~\ref{prop:single.polyhedron}(2) there exists an integral
$G$-affine pointed polyhedron $P$ such that 
$\bar{|\cP|}\subseteq\bar P$ and $\rho(P)\in\Delta$.
The desired statement now follows
from the last part of Proposition~\ref{prop:properness} 
(with $\cS = \sM(K)$) since $\trop(\cY_1)\subseteq\bar P$.
\qed

We are now in a position to state some simpler variants of 
Theorem~\ref{thm:mainthm}. However, to avoid redundancy we give the 
statements only in the strictly more general setting of multiple 
intersections.

\paragraph[Multiple intersections]
Suppose $Y$ is a smooth variety over $K$, and 
$X_1,\dots,X_m \subseteq Y$ are closed subschemes of pure codimensions 
$c_1,\dots,c_m$, with $\sum_i c_i = \dim Y$. Let $x$ be an isolated point of
$X_1 \cap \dots \cap X_m$. The \emph{local intersection number} of 
the $X_i$ at $x$ is defined to be
\[ i_K\big(x,\, X_1 \cdots X_m;\, Y\big) \coloneq 
i_K\big(D_{Y,m}(x),\, D_{Y,m}(Y) \cdot (X_1 \times \cdots \times X_m);\, Y^m\big),\]
where $D_{Y,m}:Y \to Y^m$ denotes the $m$-fold diagonal. 

Now suppose $Y=\T$, and $v$ is a (not necessarily isolated) point of 
$\Trop(X_1) \cap \dots \cap \Trop(X_m)$. Then we similarly define the 
\emph{local tropical intersection multiplicity} of the $\Trop(X_i)$ at $v$
to be
\[ i\big(v, \,\Trop(X_1) \cdots \Trop(X_m)\big) \coloneq 
i\big(D_{N_\R,m}(v), \,
D_{N_\R,m}(N_\R) \cdot (\Trop(X_1) \times \cdots \times \Trop(X_m))\big),\]
where $D_{N_\R,m}:N_\R \to (N_\R)^m$ again denotes the $m$-fold diagonal. 

\begin{rem} One can give equivalent definitions of multiple intersection
numbers inductively, if one defines multiplicities of intersections along
components using length at the generic point rather than dimension over $K$.
One then has to multiply by the degree of the residue field extension to
obtain our intersection numbers. The same equivalence holds also for
tropical intersections; see for instance 
\cite[\S~5.2]{osserman_payne:lifting}. In particular, for $m=2$ the above
definitions coincide with the definitions we have already given. This is
classical on the algebraic side, while for the tropical side one may
reduce to the algebraic side by passing to the stars at the point in
question, and expressing the tropical intersection multiplicities as
algebraic intersection multiplicities in suitable toric varieties.
\end{rem}

We have the following generalization of Theorem~\ref{thm:mainthm}:

\begin{thm}\label{thm:main-multiple} 
  Let $X_1,\dots,X_m$  be pure-dimensional closed subschemes of $\T$ with
  $\sum_i \codim(X_i)=\dim(\T)$ and $m\geq2$.  Choose polyhedral complex structures on
  the $\Trop(X_i)$.  Let $C$ be a connected component of
  $\bigcap_i \Trop(X_i)$ and suppose that $\Delta$ is an integral
  compactifying fan for $C$ such that $X(\Delta)$ is smooth.
  Let $\bar X_i$ be the closure of $X_i$
  in $X(\Delta)$ for each $i$, and let $\bar C$ be the closure of $C$ in
  $N_\R(\Delta)$.  If there are only finitely many points 
  $x\in|\bigcap_i \bar X_i|$ with $\trop(x)\in\bar C$ then
  \[ \sum_{\substack{x\in|\bigcap_i \bar X_i|\\\trop(x)\in\bar C}}
  i_K\big(x,\, \bar X_1\cdots\bar X_m;\, X(\Delta)\big) =
  \sum_{v\in C} i\big(v,\, \Trop(X_1)\cdots\Trop(X_m)\big). \]
  Furthermore, if there exists $\sigma \in \Delta$ such that $\rho(P)$ is a
  face of $\sigma$ for every polyhedron $P$ of $C$, then there are
  automatically only finitely many points $x\in|\bigcap_i \bar X_i |$ with
  $\trop(x)\in\bar C$.
\end{thm}

\pf In this proof we closely follow the statement of
Theorem~\ref{thm:mainthm}, matching our construction with its hypotheses.
The schemes $D_{\T,m}(\T)$ and $\prod_i X_i$ are pure-dimensional
closed subschemes of $\T^m$ of complementary codimension.  We have
$\Trop(D_{\T,m}(\T)) = D_{N_\R,m}(N_\R)$, which is a single polyhedron,
and $\Trop(\prod_i X_i) = \prod_i \Trop(X_i)$, which has a
polyhedral complex structure induced by the polyhedral complex structures
on the $\Trop(X_i)$.  Clearly $D_{N_\R,m}(C)$ is a connected
component of 
$D_{N_\R,m}(N_\R) \cap \prod_i \Trop(X_i) = D_{N_\R,m}(\bigcap_i\Trop(X_i))$.
We claim that $(\Delta^m,C^m)$ is a compactifying datum for
$D_{\T,m}(\T)$, $\prod X_i$, and
$D_{N_\R,m}(C)$. It is clear that 
\[C^m \cap D_{N_\R,m}(N_\R) \cap 
(\Trop(X_1) \times \cdots \times \Trop(X_m)) = 
D_{N_\R,m}(C),\]
while the fact that recession cones commute with products immediately 
implies that $\Delta^m$ is a compactifying fan for $C^m$. Finally, since
$C^m \cap \Trop(X_1) \times \cdots \times \Trop(X_m) = C^m$, we have
that $\Delta^m$ is compatible with 
$C^m \cap \Trop(X_1) \times \cdots \times \Trop(X_m)$.  Note that
$X(\Delta^m) = X(\Delta)^m$ is smooth when $X(\Delta)$ is smooth.

Since $D_{X(\Delta),m}$ is a closed immersion, the closure of 
$D_{\T,m}(\T)$ in $X(\Delta)^m$ is $D_{X(\Delta),m}(X(\Delta))$, and since
scheme-theoretic closure commutes with fiber products in this situation, the
closure of $\prod_i X_i$ is $\prod_i\bar X_i$.  Likewise, since
$D_{N_\R(\Delta),m}$ is a closed embedding, the closure of 
$D_{N_\R,m}(C)$ in $N_\R(\Delta)^m$ is $D_{N_\R(\Delta),m}(\bar C)$.
Hence there are only finitely many points of 
$|\bar{D_{\T,m}(\T)}\cap\bar{\prod_i X_i}| = 
|D_{X(\Delta),m}(\bigcap_i\bar X_i)|$ tropicalizing to 
$\bar{D_{N_\R,m}(C)} = D_{N_\R(\Delta),m}(\bar C)$.
Therefore the hypotheses of Theorem~\ref{thm:mainthm} are satisfied,
and the result follows.

Finally, if there exists $\sigma \in \Delta$ such that $\rho(P)$ is a face
of $\sigma$ for every cell $P$ of $C$, then 
$\rho(\prod_i P_i) = \prod_i \rho(P_i)$ is a face of $\sigma^m$ for every
cell $\prod_i P_i$ of $C^m$, so the
finiteness condition follows from Proposition~\ref{prop:automatic.finiteness}.
\qed

\begin{rem} \label{rem:fan.existence.3}
  By Remark~\ref{rem:cpct.fans.exist},
  Proposition~\ref{prop:basic-fan-props}(3), and the theorem on toric
  resolution of singularities, there exists an integral compactifying fan
  for $C$ such that $X(\Delta)$ is smooth.
  The condition that $\bigcap_i \bar X_i$ be finite is more subtle; it
  can certainly happen that $\bigcap_i \bar X_i$ meets the boundary
  $X(\Delta)\setminus\T$ in a positive-dimensional subset even when 
  $\bigcap_i X_i$ is finite, if the last assertion of
  Theorem~\ref{thm:main-multiple} is not applicable.
\end{rem}

\begin{rem}
  We can dispense with the smoothness hypothesis in
  Theorem~\ref{thm:main-multiple} in the case of a complete intersection.
  More precisely, suppose that each $X_i$ is the hypersurface cut out by a
  nonzero Laurent polynomial $f_i\in K[M]$.  We endow each $\Trop(X_i)$ with
  its canonical polyhedral complex structure.  Let 
  $C\subset\bigcap_{i=1}^m\Trop(X_i)$ be a connected component with its
  induced polyhedral complex structure, let $\Delta$ be an
  integral compactifying fan for $C$, and let $\bar C$ be the closure of $C$ in
  $N_\R(\Delta)$. 
  Then each $\bar X_i\cap X(\sigma)$ is again cut out by a single equation
  for $\sigma\in\Delta$ such that
  $N_\R(\sigma)$ meets $\bar C\setminus C$: indeed, for such $\sigma$ we
  have $\sigma\subset\rho(P)$ for some cell $P$ of $C$ by
  Lemma~\ref{lem:closure.P}, so we can apply~\cite[\S12]{jdr:trop_ps} to a
  cell $P'$ of $\Trop(X_i)$ containing $P$.  It  follows that  
  $\bigcap_{i=1}^m \bar X_i$ is a local complete intersection at all
  points tropicalizing to $\bar C$, and by
  Hochster's theorem the toric variety $X(\Delta)$ is Cohen-Macaulay, so
  it makes sense to define $i_K(\xi,\bar X_1\cdots\bar X_m;X(\Delta))
  = \dim_K(\sO_{\bigcap_{i=1}^m\bar X_i,\xi})$ for any isolated point
  $\xi\in\bigcap_{i=1}^m \bar X_i$ tropicalizing to $\bar C$.  In this case
  Theorem~\ref{thm:main-multiple} simply strengthens~\cite[\S12]{jdr:trop_ps}
  by adding more flexibility in the choice of fan and ground field, and is
  proved in the same way.
\end{rem}

We have the following important special case, in which no compactification
is needed:

\begin{cor} \label{cor:bounded.int}
  Let $X_1,\dots,X_m$ be pure-dimensional closed subschemes of $\T$ with
  $\sum_i \codim(X_i) =\dim(\T)$.  
  Let $C$ be a connected component of $\bigcap_i \Trop(X_i)$, and suppose
  that $C$ is bounded.  Then there are only finitely many points 
  $x\in|\bigcap_i X_i|$ with $\trop(x)\in C$, and 
  \[ \sum_{\substack{x\in|\bigcap_i X_i|\\\trop(x)\in C}}
  i_K(x,\, X_1 \cdots X_m;\, \T) =
  \sum_{v\in C} i(v,\, \Trop(X_1)\cdots\Trop(X_m)). \]
\end{cor}

\pf Apply Theorem~\ref{thm:main-multiple} with $\Delta = \{\{0\}\}$.\qed

\begin{rem}
  Suppose that $K$ is trivially-valued.  If $X\subset X(\Delta)$ is a
  closed subscheme then 
  $\trop(x)\in\Djunion_{\sigma\in\Delta} \pi_\sigma(0)\subset N_\R(\Delta)$
  for every $x\in|X|$.  In particular, the compactification required for
  Theorem~\ref{thm:main-multiple} is still necessary in this situation.
  On the other hand, Corollary~\ref{cor:bounded.int} is exactly the same
  as Theorem~\ref{thm:compat.proper} in the trivially-valued case since if
  $C$ is bounded then $C = \{0\}$.
\end{rem}

\section{An example} \label{sec:the.example}

  The following example is meant to illustrate Theorem~\ref{thm:mainthm}.
  Assume that $K$ is complete and nontrivially valued.
  Let $M = N = \Z^3$ and let $x = x^{-e_1}, y = x^{-e_2}, z = x^{-e_3}$, where
  $e_1,e_2,e_3$ is the standard basis of $\Z^3$.  We have
  $\T = \Spec(K[x^{\pm 1}, y^{\pm 1}, z^{\pm 1}])\cong\G_m^3$, and for 
  $\xi\in|\T|$ we have
  $\trop(\xi) = (\val(x(\xi)),\val(y(\xi)),\val(z(\xi)))$ according to our
  sign conventions.  

  Let $X\subset\T$ be the curve defined by the  equations 
  \begin{equation} \label{eq:nonred.cubic}
    (y-1)^2 = x(x-1)^2 \qquad (x-1)(z-1) = 0 \qquad (z-1)^2 = 0. 
  \end{equation}
  This is a slight simplification of the degeneration of a family of twisted 
  cubic curves found in~\cite[Example~III.9.8.4]{hartshorne:ag}.
  This curve has a non-reduced point at $(1,1,1)$ and is reduced
  everywhere else.  Hence $X$ is not a local complete intersection
  at $(1,1,1)$.  The tropicalization of $X$ coincides with the
  tropicalization of the underlying reduced curve $X^{\red}$, which is a nodal cubic
  curve in the $(x,y)$ plane; one computes 
  $\Trop(X) = \Trop(X^{\red})$ using the Newton polytope of the defining
  equation $(y-1)^2 = x(x-1)^2$.  The tropicalization equal to the union of the
  rays $R_1 = \R_{\geq 0}\cdot e_1$, 
  $R_2 = \R_{\geq 0}\cdot e_2$, and 
  $R_3 = \R_{\geq 0}\cdot (-2e_1-3e_2)$; these rays have tropical
  multiplicities $2,3$, and $1$, respectively.  See
  Figure~\ref{fig:nonci}.

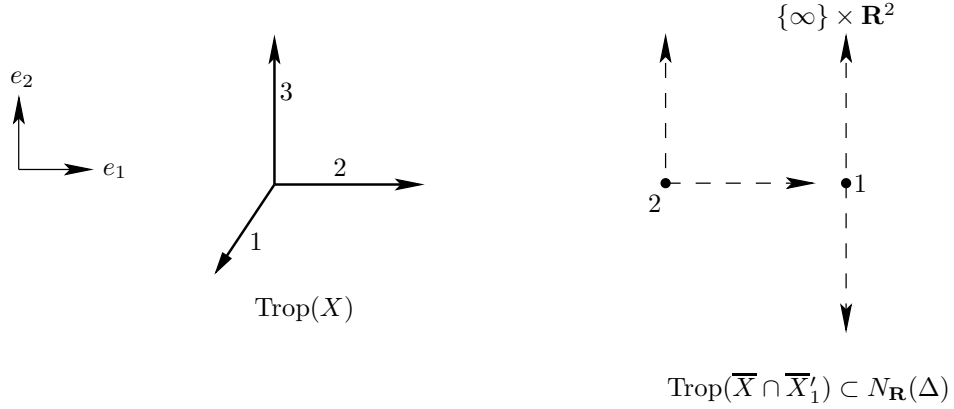
  \genericfig[ht]{nonci}{On the left: the tropicalization of the curve $X$
  defined by~\eqref{eq:nonred.cubic}, which is contained in
  $\spn(e_1,e_2)$; the $e_3$ direction is orthogonal to the page.  On
  the right is the tropicalization of $\bar X\cap\bar X\p_1$, with the
  numbers indicating the algebraic intersection multiplicities.}

  Let $X'_a\subset\T$ be the plane defined by $y = a$ for $a\in K$ with
  $\val(a) = 0$.  Then $\Trop(X'_a)$ is the plane spanned by $e_1$ and
  $e_3$, and $\Trop(X)\cap\Trop(X'_a) = R_1$.  
  The intersection of $\Trop(X'_a) + \epsilon e_2$ with
  $\Trop(X)$ is the point $\epsilon e_2$ counted with multiplicity $3$;
  hence the stable tropical intersection $\Trop(X)\cdot\Trop(X'_a)$ is the
  point $0$ counted with multiplicity $3$.

  Let $\Delta = \{\{0\}, R_1\}$.  This is a compactifying fan for $R_1$.
  We have $N_\R(\Delta) = N_\R\djunion(N_\R/\spn(e_1))$ and
  $X(\Delta)\cong\Spec(K[x,y^{\pm 1},z^{\pm 1}])\cong\A^1\times\G_m^2$; if
  we identify 
  $N_\R/\spn(e_1)$ with $\{\infty\}\times\R^2$ then the tropicalization
  map $\trop:|X(\Delta)|\to N_\R(\Delta)$ again can be written
  $\trop(\xi) = (\val(x(\xi)),\val(y(\xi)),\val(z(\xi)))$, since
  $\val(0) = \infty$.  The closure $\bar R_1$ of $R_1$ in $N_\R(\Delta)$ is 
  $R_1\djunion\{(\infty,0,0)\}$, the closure $\bar X$ of $X$ in
  $X(\Delta)$ is also given by~\eqref{eq:nonred.cubic}, and the closure of
  $\bar X\p_a$ of $X_a'$ in $X(\Delta)$ is also given by $\{y=a\}$.

  Let us calculate $\bar X\cdot\bar X\p_a$.  The scheme-theoretic
  intersection $\bar X\cap\bar X\p_a$ is defined by 
  the ideal
  \[ I_a = \big(x(x-1)^2 - (y-1)^2,~ (x-1)(z-1),~ (z-1)^2,~ y-a\big); \]
  hence $\bar X\cap\bar X\p_a$ is supported on the points of the form 
  $(r, a, 1)\in|\T|$, where $r$ is a root of the cubic polynomial 
  $q_a(x) = x^3 - 2x^2 + x - (a-1)^2$.

  Suppose first that $a = 1$, so $q_1(x) = x(x-1)^2$ and 
  $\bar X\cap\bar X\p_1$ is supported on the points
  $\xi_1 = (1,1,1)$ and $\xi_0 = (0,1,1)$.  The point
  $\xi_0$ is reduced in $\bar X\cap\bar X\p_1$ and is a smooth
  point of both $\bar X$ and $\bar X\p_1$; hence 
  $i_K(\xi_0,\bar X\cdot\bar X\p_1;X(\Delta)) = 1$.  
  We identify the completed local ring of $\T$ at $\xi_1$ with
  $B\coloneq K\ps{x_1,y_1,z_1}$, where $x_1 = x + 1$, $y_1 = y + 1$, and 
  $z_1 = z + 1$.  Then the completed local ring of $X$ at $\xi_1$ is 
  \[ A = K\ps{x_1,y_1,z_1}/(x_1^2(x_1+1),\, x_1z_1,\, z_1^2)
  \cong K\ps{x_1,y_1,z_1}/(x_1^2,\, x_1z_1,\, z_1^2), \]
  and the completed local ring of $X'_1$ at $\xi_1$ is 
  $A' = K\ps{x_1,z_1}$.  Hence the local ring of $X\cap X'_1$ at $\xi_1$
  is  
  \[ A\tensor_B A' \cong 
  K\ps{x_1,y_1,z_1}/(x_1^2,\, x_1z_1,\, z_1^2,\,y_1)
  \cong K\ps{x_1,z_1}/(x_1,z_1)^2, \]
  which is an Artin ring of dimension $3$ over $K$.
  We have a resolution 
  \[ 0 \To B \overset{\cdot y_1}\To B\To A' \To 0, \]
  so the groups $\Tor_i^B(A,A')$ are calculated by the complex
  \[ 0 \To A \overset{\cdot y_1}\To A \To 0. \]
  Hence $\Tor_i^B(A,A') = 0$ for $i > 1$, and 
  $\Tor_1^B(A,A')$ is identified with the space of $y_1$-torsion in $A$.
  It is not hard to see that $\Tor_1^B(A,A')$ is spanned over $K$ by
  $y_1z_1$, so
  \[ i_K(\xi_1,\,\bar X\cdot\bar X\p_1,\,X(\Delta)) 
  = \dim_K(A\tensor_B A') - \dim_K(\Tor_1^B(A,A')) 
  = 3 - 1 = 2. \]
  Therefore 
  \[ 3 = \sum_{v\in\bar R_1} i(v,\, \Trop(X)\cdot\Trop(X_1')) 
  = \sum_{\xi\in|\bar X\cap\bar X\p_1|} i_K(\xi,\,\bar X\cdot\bar X\p_1,\,X(\Delta))
  = 2 + 1, \]
  as in Theorem~\ref{thm:mainthm}.  Note that we would have gotten the
  wrong number on the right side of the above equation if we had
  na\"ively defined the intersection number at $\xi_1$ as the
  dimension of the local ring of $X\cap X_1'$, or if we had not passed to
  the toric variety 
  $X(\Delta)$ which compactifies the situation in the direction of $R_1$.%
  \footnote{Coincidentally, if we had done neither of these things then
    the intersection numbers would coincide in this example.}

  Now suppose that $a\neq 1$ (but still $\val(a)=0$).  In this case 
  $\bar X\cap\bar X\p_a = X\cap X_a'$, and every point 
  $\xi\in|X\cap X_a'|$ is a smooth point of both $X_a$ and $X_a'$,
  so $i_K(\xi,\bar X\cdot\bar X\p_a; X(\Delta))$ is equal to the dimension of
  the local ring of $X\cap X'_a$ at $\xi$.  Writing $\xi = (r,a,1)$, we have
  $\trop(\xi) = (\val(r),0,0)$.  The possible values for $\val(r)$ are
  easily calculated from the Newton polygon of 
  $q_a(x) = x^3 - 2x^2 + x - (a-1)^2$; the result of this calculation is
  that there are two points $\xi$ (counted with multiplicity) with
  $\trop(\xi) = (0,0,0)$, and one with $\trop(\xi) = (\val(a-1),0,0)$.  In
  particular, $\trop(\xi)$ can lie \emph{anywhere} on $\bar R_1\cap N_G$, so we
  cannot strengthen Theorem~\ref{thm:mainthm} in such a way as to pinpoint
  $\Trop(\bar X\cap\bar X\p,\Delta)$ more precisely.

\bibliographystyle{thesis}
\bibliography{int_nums}

\providecommand{\noopsort}[1]{}\def\cprime{$'$}
\providecommand{\bysame}{\leavevmode\hbox to3em{\hrulefill}\thinspace}
\providecommand{\MR}{\relax\ifhmode\unskip\space\fi MR }
% \MRhref is called by the amsart/book/proc definition of \MR.
\providecommand{\MRhref}[2]{%
  \href{http://www.ams.org/mathscinet-getitem?mr=#1}{#2}
}
\providecommand{\href}[2]{#2}
\begin{thebibliography}{EGAIII${}_1$}

\bibitem[Ber90]{berkovich:analytic_geometry}
V.~G. Berkovich, \emph{Spectral theory and analytic geometry over
  non-{A}rchimedean fields}, Mathematical Surveys and Monographs, vol.~33,
  American Mathematical Society, Providence, RI, 1990.

\bibitem[Ber93]{berkovich:etalecohomology}
\bysame, \emph{{\'E}tale cohomology for non-{Archimedean} analytic spaces},
  Inst. Hautes {\'E}tudes Sci. Publ. Math. \textbf{78} (1993), 5--161.

\bibitem[BGR84]{bgr:nonarch}
S.~Bosch, U.~G{\"u}ntzer, and R.~Remmert, \emph{Non-{A}rchimedean analysis},
  Grundlehren der Mathematischen Wissenschaften [Fundamental Principles of
  Mathematical Sciences], vol. 261, Springer-Verlag, Berlin, 1984.

\bibitem[BPR11]{bpr:nonarch_trop}
M.~Baker, S.~Payne, and J.~Rabinoff, \emph{Non-{A}rchimedean geometry,
  tropicalization, and metrics on curves}, 2011, Preprint available at
  \url{http://arxiv.org/abs/1104.0320}.

\bibitem[Con99]{conrad:irredcomps}
B.~Conrad, \emph{Irreducible components of rigid spaces}, Ann. Inst. Fourier
  (Grenoble) \textbf{49} (1999), no.~2, 473--541.

\bibitem[Duc07]{ducros:relative_dimension}
A.~Ducros, \emph{Variation de la dimension relative en g\'eom\'etrie analytique
  {$p$}-adique}, Compos. Math. \textbf{143} (2007), no.~6, 1511--1532.

\bibitem[Duc11]{ducros:flatness}
\bysame, \emph{Flatness in non-archimedean analytic geometry}, 2011, Preprint
  available at \url{http://arxiv.org/abs/1007.4259}.

\bibitem[EGAIII${}_1$]{egaiii_1}
A.~{\noopsort{EGAIII_1}}Grothendieck, \emph{\'{E}l\'ements de g\'eom\'etrie
  alg\'ebrique. {III}. \'{E}tude cohomologique des faisceaux coh\'erents. {I}},
  Inst. Hautes \'Etudes Sci. Publ. Math. (1961), no.~11, 167.

\bibitem[FS97]{fulton_sturmfels:itheory}
W.~Fulton and B.~Sturmfels, \emph{Intersection theory on toric varieties},
  Topology \textbf{36} (1997), no.~2, 335--353.

\bibitem[Har77]{hartshorne:ag}
R.~Hartshorne, \emph{Algebraic geometry}, Springer-Verlag, New York, 1977,
  Graduate Texts in Mathematics, No. 52.

\bibitem[MS09]{maclagan_sturmfels:book}
D.~Maclagan and B.~Sturmfels, \emph{Introduction to tropical geometry}, 2009,
  Preprint available at
  \url{www.warwick.ac.uk/staff/D.Maclagan/papers/TropicalBook.pdf}.

\bibitem[OP10]{osserman_payne:lifting}
B.~Osserman and S.~Payne, \emph{Lifting tropical intersections}, 2010, Preprint
  available at \url{http://arxiv.org/abs/1007.1314}.

\bibitem[Pay09]{payne:analytification}
S.~Payne, \emph{Analytification is the limit of all tropicalizations}, Math.
  Res. Lett. \textbf{16} (2009), no.~3, 543--556.

\bibitem[Rab10]{jdr:trop_ps}
J.~Rabinoff, \emph{Tropical analytic geometry, {N}ewton polygons, and tropical
  intersections}, 2010, Preprint available at
  \url{http://arxiv.org/abs/1007.2665}.

\bibitem[Roh11]{rohrer:completions_fans}
F.~Rohrer, \emph{Completions of fans}, 2011, Preprint available at
  \url{http://arxiv.org/abs/1107.2483}.

\end{thebibliography}
\bigskip~\bigskip

\end{document}